\documentclass[11pt]{article}
\usepackage[psamsfonts]{amsfonts}
\usepackage{amsmath,amssymb,amsthm}
\usepackage[dvips]{graphicx}
\usepackage{eepic}
\newtheorem{THM}{Theorem}[section]

\newtheorem{REM}[THM]{Remark}

\newtheorem{PRP}[THM]{Proposition}

\newtheorem*{ASSLLN}{Assumption (LLN)}
\newtheorem*{ASSCLT}{Assumption (CLT)}

\def\UseSection{
        \numberwithin{equation}{section}
	\theoremstyle{plain}
        \newtheorem{theorem}    {Theorem}[section]
        \DefineTheorems 
}

\def\DefineTheorems{
	
	\newtheorem{lemma}      [theorem] {Lemma}
	
	\newtheorem{prop}       [theorem] {Proposition}
	
	\newtheorem{cor}        [theorem] {Corollary}

	\theoremstyle{definition}
	\newtheorem{defn}       [theorem] {Definition}

	\theoremstyle{definition}

}

\newcommand{\bt}   {\begin{theorem}}
\newcommand{\et}   {\end  {theorem}}
\newcommand{\bl}   {\begin{lemma}}
\newcommand{\el}   {\end  {lemma}}
\newcommand{\bp}   {\begin{prop}}
\newcommand{\ep}   {\end  {prop}}
\newcommand{\bc}   {\begin{cor}}
\newcommand{\ec}   {\end  {cor}}
\newcommand{\bd}   {\begin{defn}}
\newcommand{\ed}   {\end  {defn}}

\newcommand{\ba}   {\begin{array}}
\newcommand{\ea}   {\end  {array}}
\newcommand{\be}   {\begin{enumerate}}
\newcommand{\ee}   {\end  {enumerate}}
\newcommand{\bi}   {\begin{itemize}}
\newcommand{\ei}   {\end  {itemize}}

\def\eq#1\en{\begin{equation}#1\end{equation}}  
\def\eqsplit#1\ensplit{
	\begin{equation}\begin{split}#1\end{split}\end{equation}
	}
\def\eqalign#1\enalign{
	\begin{align}#1\end{align}
	}
\def\eqmul#1\enmul{
	\begin{multline}#1\end{multline}
	}
\newcommand{\eqarrstar} {\begin{eqnarray*}} 
\newcommand{\enarrstar} {\end{eqnarray*}} 
\newcommand{\eqarray}   {\begin{eqnarray}} 
\newcommand{\enarray}   {\end{eqnarray}}

\newcommand{\lbeq}[1]  {\label{e:#1}}
\newcommand{\refeq}[1] {\eqref{e:#1}}    

\makeatletter
\newcommand{\labelcounter}[2]{{%
	\stepcounter{#1}
	\protected@write\@auxout{}%
	{\string\newlabel{#2}{{\csname the#1\endcsname}{\thepage}}}%
	{\ref{#2}}
	}}
\makeatother
\newcommand{\sss}   { \scriptscriptstyle } 

\newcommand{\Cbold} {{\mathbb C}}  
  
\newcommand{\Ebold} {{\mathbb E}}

\newcommand{\Nbold} {{\mathbb N}}

\newcommand{\Pbold} {{\mathbb P}}
\newcommand{\Qbold} {{\mathbb Q}}
\newcommand{\Rbold} {{\mathbb R}}

\newcommand{\Zbold} {{\mathbb Z}}

\newcommand{\Zd}    {{ {\Zbold}^d }}

\newcommand{\spose}[1] {{\hbox to 0pt{#1\hss}} }
\newcommand{\ltapprox} {\mathrel{\spose{\lower 3pt\hbox{$\mathchar"218$}}
 \raise 2.0pt\hbox{$\mathchar"13C$}}}
\newcommand{\gtapprox} {\mathrel{\spose{\lower 3pt\hbox{$\mathchar"218$}}
 \raise 2.0pt\hbox{$\mathchar"13E$}}}

\UseSection  
\newcommand{\extend}[1]{#1}
\newcommand{\short}[1]{}

\usepackage[usenames]{color}

\newcommand{\ch}[1]{#1}
\newcommand{\ci}[1]{#1}
\newcommand{\mQ}{\mathbb{Q}}

\newcommand{\hlf}{\frac{1}{2}}
\newcommand{\ra}{\rightarrow}

\newcommand{\prob}{{\mathbb P}}

\newcommand{\shift}   {\!\!\!\!}

\newcounter{countC}  
\newcounter{countR}  
\newcommand{\R}{\Rbold}
\newcommand{\Z}{\Zbold}
\newcommand{\N}{\Nbold}
\newcommand{\C}{\Cbold}

\newcommand{\conn}{\longrightarrow}

\newcommand{\nn}{\nonumber}

\newcommand{\smallsup}[1] {{\scriptscriptstyle{({#1}})}}

\newcommand{\walk}{\vec{\omega}}
\newcommand{\ewalk}{\vec{\eta}}
\newcommand{\walkvec}[2]{\vec{\omega}^{\smallsup{#1}}_{#2}}
\newcommand{\walkcoor}[2]{\omega^{\smallsup{#1}}_{#2}}

\newcommand{\zerospeed}{\theta_{\sss \varnothing}}
\newcommand{\zerovar}{\Sigma_{\sss \varnothing}}
\newcommand{\zerop}{p^{\sss 0}}

\newcommand{\floor}[1]{\lfloor #1 \rfloor}

\newcommand{\Var}{{\rm Var}}
\newcommand {\convd}{\stackrel{d}{\Longrightarrow}}
\newcommand {\convp}[1]{~\overset{#1}{\conn}~}

\newcommand{\mc}[1]{\mathcal{#1}}
\newcommand{\mE}{\mathbb{E}}
\newcommand{\mP}{\mathbb{P}}

\newcommand{\eqn}[1]{\eq #1 \en}
\newcommand{\vep}{\varepsilon}
\newcommand{\veta}{\vep_\beta}
\newcommand{\sN}{{\sss N}}

\title  {
        An expansion for self-interacting random walks\extend{:\\
        extended version}
        }

\author{
Remco van der Hofstad\footnote{Department of Mathematics and
Computer Science, Eindhoven University of Technology, P.O.\ Box
513, 5600 MB Eindhoven, The Netherlands. E-mail {\tt
rhofstad@win.tue.nl}}
\and
Mark Holmes\footnote{Department of Statistics,
The University of Auckland, Private Bag 92019, Auckland 1142,
New Zealand. E-mail {\tt mholmes@stat.auckland.ac.nz}}
}

\begin{document}

\maketitle

    \begin{abstract}
    We derive a perturbation expansion for general self-interacting
    random walks, where steps are made on the basis of the
    history of the path.  Examples of models where this
    expansion applies are reinforced random walk, excited random walk,
    the true (weakly) self-avoiding walk, loop-erased random walk,
    and annealed random walk in random environment.
    In this paper we show that the expansion gives rise to useful
    formulae for the speed and variance of the random walk, when
    these quantities are known to exist.  The results and formulae
    of this paper have been used elsewhere by the authors to prove
    monotonicity properties for the speed (in high dimensions) of
    excited random walk and related models, and certain models
    of random walk in random environment.
    We also derive a law of large numbers and central limit theorem
    (with explicit error terms) directly from this expansion,
    under strong assumptions on the expansion coefficients.
    The assumptions are shown to be satisfied by excited random
    walk in high dimensions with small excitation parameter,
    a model of reinforced random walk with underlying drift
    and small reinforcement parameter, \ch{and certain models
    of random walk in random environment under strong ellipticity
    conditions.}
    \short{Certain technical proofs are deferred to the
    extended version of our paper \cite{HofHol10b}.}
    \extend{This is the extended version of the paper
    \cite{HofHol10a}, where we provide all proofs.}
\end{abstract}


\section{Introduction}
\label{sec-intro}
Recently, many models of random walks with a
certain self-interaction have been introduced. A few examples are self-reinforced
random walks \cite{DurKesLim02,Pema88, Roll02}, excited random walks
\cite{BenWil03, Kozm03, Kozm05, Zern05, Zern06},
true-self avoiding walks and loop-erased random walks.
Proofs in these models often rely on martingale methods,
or explicit comparisons to random walk properties.
In some of the examples, laws of large numbers are
derived. The difficulty is that the limiting parameters
are rather implicit, so that it is hard to derive analytical
properties of them. For example, it is quite reasonable to
assume that the drift for excited random walk is monotone increasing
in the excitement parameter for each $d\ge2$, but a proof of this fact is currently
missing. Similarly, it has not been proved that the speed for once-reinforced
random walk on the tree is monotone decreasing in the reinforcement
parameter (see \cite{DurKesLim02}).
See \cite{Pema07} for a survey of self-interacting random walks with
reinforcement.

{\extend{In the past decades, the lace expansion has proved to be an
extremely useful technique to investigate a variety of models
above their upper-critical dimension,
where Gaussian limits are expected.
Examples are self-avoiding walks above 4 dimensions \cite{BrySpe85, HarSla92b, Slad87, Slad88, Slad89},
lattice trees above 8 dimensions \cite{DerSla97, DerSla98,HarSla90b, Holm08}, the contact process
above 4 dimensions \cite{HofSak04, HofSak08}, oriented percolation above 4
dimensions \cite{HofSla03b, NguYan93, NguYan95}, and percolation above 6
dimensions \cite{HarSla90a, HarSla00a, HarSla00b}.  An essential ingredient in the
proofs is the fact that the above models are {\it self-repellent}. There are
many more models where a Gaussian limit is expected above a certain upper
critical dimension, but using the lace expansion for these models is hard
as they are not strictly self-repellent. In this paper, we perform a first
step for a successful application of the lace expansion methodology, namely,
we derive the expansion for general self-interacting random walks.
The goal is to use this expansion for some of the simpler self-interacting stochastic processes
available.}}

We will study a particular version of once-reinforced random walk, where the
initial weights are such that the corresponding random walk has a {\it non-zero drift}.  A similar situation was investigated in \cite{DurKesLim02}, where once-reinforced
random walk was investigated on the tree. We expect that our method can be adapted to the tree setting to reprove
some results in \cite{DurKesLim02}, when the reinforcement parameter is sufficiently small.  We also study {\it excited random walk}, where the random walker has a drift
in the direction of the first component each time when the walker visits a
new site.  It was shown that this process has ballistic behaviour when $d\geq 2$ in
\cite{BenWil03,Kozm03, Kozm05}, while there is no ballistic behaviour in one dimension.  For a third application of our method we study random walk in (partially) random environments, similar to those considered in \cite{BSZ03}.


In this paper we give a (self-contained) proof of a law of large numbers ($d\ge 6$) and a central limit theorem $d\ge 9$) when the excitation parameter is sufficiently small.  We also derive a law of large numbers and central limit theorem for the once-reinforced random walk with drift when the reinforcement is sufficiently small compared to the drift.  These results were completed in 2006.  Since then, substantial progress has been made on these two models.  Using renewal techniques, a strong law of large numbers and invariance principle has been proved for the excited random walk in dimensions $d\ge2$ \cite{BR07}, while laws of large numbers and local central limit theorems are obtained for a large class of ballistic self-interacting random walks (including the reinforced random walk with drift in all dimensions) in \cite{IV08}.  We prove similar results for annealed random walks in partially random environment similar to those in \cite{BSZ03} in the perturbative regime, with the difference being that the probability of taking a step in each coordinate may be random.  We believe that our results for random walk in random environment are new.

{\extend{The renewal techniques often give strong results (as described above), but currently do not provide much insight into how the results depend on the underlying parameters.  As is done in this paper, our expansion can be used independently to prove (sometimes weaker) results in the perturbative regime.  In doing so we obtain formulae and estimates of error terms for some of the relevant quantities of interest.  This is one of the main advantages of our method, but as illustrated by recent applications (see Section \ref{sec-applic}) we see a {\em combination} of our expansion with renewal and ergodic methods to be highly informative.
}}

\section{The main results}
\label{sec-results}
We start by introducing some notation. A path $\omega$ is a sequence
$\{\omega_i\}_{i=0}^\infty$ for which $\omega_i\in \Z^d$
for all $i\geq 0$.  We obtain {\it random walk} when the random vector
$\{\omega_{i+1}-\omega_i\}_{i=0}^{\infty}$ is an i.i.d.\ sequence.
We let $\Pbold$ be the law of a random walk law starting at
the origin. We write $\walk_n$ for the vector
    \eq
    \walk_n=(\omega_0, \ldots, \omega_n),
    \en
that is, for the first $n$ positions of the walk and its starting point. Let
    \eq
    D(x)=\Pbold(\omega_1=x)
    \en
be the random walk transition probability, so that
    \eq
    \lbeq{SRWprod}
    \Pbold(\walk_n=(x_0, x_1, \ldots, x_n))
    =\prod_{i=0}^{n-1} D(x_{i+1}-x_{i}).
    \en
    We restrict our attention to $D$ with finite range $L$ so that $\sum_{x:|x|>L}D(x)=0$ and all moments of $D$ exist.
For self-interacting random walks, a similar expression to \refeq{SRWprod} is valid, but
the term appearing in the product may depend on the history of the
path. \ch{
Let $\Qbold_{(x_0)}$ denote the law of a self-interacting random walk $\walk=(\omega_0,\omega_1,\dots)$ started at $\omega_0=x_0$, i.e.,
    \eq
    \lbeq{SIRP}
    \Qbold_{(x_0)}(\walk_n=(x_0, x_1, \ldots, x_n))
    =\prod_{i=0}^{n-1} p^{\vec{x}_i}(x_{i},x_{i+1}),
    \en
where
\[p^{\vec{x}_i}(x_i,x_{i+1})=\Qbold_{(x_0)}(\omega_{i+1}=x_{i+1}|\walk_i=\vec{x}_i).\]
In other words, for a general path $\vec{x}_i$, we write $p^{\vec{x}_i}(x_i, x_{i+1})$
for the conditional probability that the walk steps from $x_i$ to
$x_{i+1}$, given the history of the entire path $\vec{x}_i=(x_0, \ldots, x_i)$.
It is crucial to our analysis that our self-interacting random walk law is translation invariant, i.e.~for all $\vec{x}_n$,
\[\Qbold_{(x_0)}(\walk_n=(x_0, x_1, \ldots, x_n))= \Qbold_{(o)}(\walk_n=(o, x_1-x_0, \ldots, x_n-x_0)) .\]
We henceforth write $\Qbold=\Qbold_{(o)}$ and drop the dependence on the starting point $x_0$ from the notation when the history $\vec{x}_n$ of the path is given, e.g.~$\Qbold(\cdot|\walk_n=\vec{x}_n)=\Qbold_{(x_0)}(\cdot|\walk_n=\vec{x}_n)$.}

The goal of this paper is to investigate the {\it two-point function}
    \eq
    \lbeq{twopoint}
    c_n(x) =\Qbold(\omega_n=x).
    \en


In this paper, we will derive an expansion for the two-point function in full generality.  However, for the analytical results we will focus on directed once-edge-reinforced random walks, excited random walks, and random walks in partially random environments. In Sections \ref{sec-OERRW}, \ref{sec-ERW}, and \ref{sec-RWRE}
below, we will define the models and state the results.

\subsection{Once edge-reinforced random walk with drift}
\label{sec-OERRW}
In this section, we introduce an example
of a once edge-reinforced random walk with drift.
For a directed edge $b$, denote the number of times the edge $b$ is traversed up to time $t$ by
\eq
    \ell_t(b)=\sum_{i=1}^t I_{\{(\omega_{i-1},\omega_i)=b\}},
    \en
where $I_{A}$ denotes the indicator of the event $A$,
and let $t\mapsto \beta_t$ be a sequence of $\R$-valued reinforcement parameters.  We use $w_s(b)$ to denote the {\it weight of the
edge} $b$ at time $s$.  The main assumption for our reinforced random walk is that
$w_0(b)$ is translation invariant, and that
    \eq
    \lbeq{mainass}
    \sum_x xw_0(0,x)\neq 0.
    \en
Define $w_s(b)$ recursively by
    \eq
    \lbeq{weights}
    w_t(b)= w_{t-1}(b)+I_{\{(\omega_{t-1},\omega_t)=b\}} \beta_{\ell_t(b)}.
    \en
We define a directed version of {\ci edge-reinforced random walk} (ERRW) by setting
    \eq
    p^{\walk_i}(x_{i},x_{i+1})=\frac{w_i(x_i, x_{i+1})}{\sum_{y} w_i(x_i, y)}.
    \en

We will deal with directed {\it once-reinforced random walks}, where $\beta_{t}=\beta\delta_{t,1}$ is taken sufficiently small, however our results extend to directed {\it boundedly-reinforced random walks},
where we assume that
    \eq
    \lbeq{betadef}
    \beta=\sum_{t=0}^{\infty} |\beta_t|<\infty, \quad \text{is sufficiently small}.
    \en
The parameters $\beta_s$ are allowed to be negative (provided $w_0(b)+\sum_{t=1}^m\beta_t$ remains bounded away from $0$).
  Note that \refeq{mainass} implies that the random walk distribution arising for $\beta=0$ has {\it non-zero drift}.


We denote by $\Qbold_\beta$ the distribution of the above
once-reinforced random walk with drift, and we let $\Ebold_\beta$
denote expectation with respect to $\Qbold_\beta$. We denote by
$\Var_\beta(\omega_n)$ the covariance matrix of the random vector
$\omega_n$ under the measure $\Qbold_\beta$.  We also denote convergence in distribution
by $\convd$, convergence in probability under the law $\prob$ by $\convp{\prob}$ and write ${\cal N}(0,\Sigma)$ for the multivariate normal
distribution with mean the zero vector and covariance matrix
$\Sigma$.

\begin{theorem}[A CLT for finitely reinforced random walk with drift]
\label{thm-1OERRW} Fix $d\geq 1$ and assume
\refeq{mainass}. There exist $\beta_0 = \beta_0(d, w_0)>0$,
$\theta=\theta(\beta, w_0,d)\in [-1,1]^d$ and finite $\Sigma=
\Sigma(\beta, w_0,d)$ such that,
for all $\beta \leq \beta_0$,\\
(a)
\eq
\lbeq{mean}
    \Ebold_\beta[\omega_n]
    =\theta n [1+O(\frac{1}{n})].
\en
(b)
\eq
\lbeq{var}
    \Var_\beta(\omega_n)
    = \Sigma n [1+O(\frac{1}{n})].
\en
(c) $\omega_n$ satisfies a central limit theorem under
$\Qbold_\beta$, that is,
\eq
\lbeq{CLT}
    \frac{\omega_n-\theta n}{\sqrt{n}} \convd {\cal N}(0,\Sigma).
\en
%
\end{theorem}

As noted in the introduction, this result has
since been strengthened in \cite{IV08}, although without error estimates.
Parts of our methods apply to once-reinforced random
walk where the initial weights induce no drift. However, we are currently
unable to prove the bounds on the expansion coefficients, one of the crucial
steps in the analysis. {\extend{We shall comment on this issue in more detail in Section
\ref{sec-disc-bds} below.}} In \ch{Section \ref{sec-expspeed}-\ref{sec-expvar},}
we shall further give formulas for the speed and variance appearing in Theorem \ref{thm-1OERRW}.

%

\subsection{Excited random walk}
\label{sec-ERW}
In this section, we introduce excited random walk (ERW), which is
the second model to which we shall apply our expansion method.
It is defined for $\beta \in [0,1]$ by taking
    \eq
    p^{\walk_i}(x_{i},x_{i+1})=p_0(x_{i+1}-x_i)
    \ch{I_{\{x_i \in \walk_{i-1}\}}}
    +p_\beta (x_{i+1}-x_i)\big[1-I_{\{x_i \in \walk_{i-1}\}} \big],
    \en
where $\{x_i \in \walk_{i-1}\}$ denotes the event that $x_i=\omega_j$ for some
$0\leq j\leq i-1$, and where
    \eq
    p_0(x)=\frac{1}{2d} I_{\{|x|=1\}}
    \en
is the nearest-neighbour step distribution and
    \eq
    p_\beta(x) =\frac{1+\beta e_1\cdot x}{2d} I_{\{|x|=1\}}.
    \en
Here $e_1=(1,0, \ldots, 0)$ and $x\cdot y$ is the inner-product between $x$ and $y$.
In words, the random walker gets excited and has a positive drift in the direction of the first
coordinate each time he/she visits a new site.

That ERW has a positive drift (in the sense of a lower bound) was established for ERW in $d\geq 4$
in \cite{BenWil03}, for $d=3$ in \cite{Kozm03}, and for $d=2$ in
\cite{Kozm05}.  For $d=1$, it is known that ERW is recurrent and diffusive (except the trivial case $\beta=1$) \cite{Davi99}.  Many generalisations of this model, described in terms of cookies, have also been studied (see for example \cite{Zern05}, \cite{AntRed05}, \cite{BS08}).

We denote by $\Qbold_\beta$ the distribution of the above
excited random walk started at the origin, and we let $\Ebold_\beta$
denote expectation with respect to $\Qbold_\beta$. We denote by
$\Var_\beta(\omega_n)$ the covariance matrix of the random vector
$\omega_n$ under the measure $\Qbold_\beta$.

Our main result for excited random walk is the following theorem:

\begin{theorem}[A CLT for ERW above 8 dimensions]
\label{thm-1ERW}
Fix $d>8$. Then, there exists $\beta_0 = \beta_0(d)>0$,
$\theta=(\theta_1(\beta,d), 0,\ldots, 0)$ and finite
$\Sigma=\Sigma(\beta, d)$ such that,
for all $\beta \leq \beta_0$,\\
(a)
\eq
\lbeq{meanERW}
    \Ebold_\beta[\omega_n]
    =\theta n [1+O(\frac{1}{n})].
\en
(b)
\eq
\lbeq{varERW}
    \Var_\beta(\omega_n)
    = \Sigma n [1+O(\frac{\log n}{n^{1 \wedge \frac{d-7}{2}}})].
\en
(c) $\omega_n$ satisfies a central limit theorem under
$\Qbold_\beta$, that is,
\eq
\lbeq{CLTERW}
    \frac{\omega_n-\theta n}{\sqrt{n}} \convd {\cal N}(0,\Sigma).
\en
%
\end{theorem}


Unfortunately, our methods do not apply to general $d\geq 2$.
However, when $d>5$, we {\it can} prove a weak law of large numbers:

\begin{theorem}[A LLN for ERW above 5 dimensions]
\label{thm-2ERW}
Fix $d>5$. Then, there exists $\beta_0 = \beta_0(d)>0$ and
$\theta=(\theta_1(\beta,d), 0,\ldots, 0)$ such that,
for all $\beta \leq \beta_0$,\\
(a)
\eq
\lbeq{meanERWLLN}
    \Ebold_\beta[\omega_n]
    =\theta n [1+O(\frac{\log{n}}{n^{1 \wedge \frac{d-5}{2}}})].
\en
(b) $\omega_n$ satisfies a law of large numbers under
$\Qbold_\beta$, that is,
\eq
\lbeq{LLNERW}
    \frac{\omega_n}{n} \convp{\Qbold_\beta} \theta.
\en
%
\end{theorem}
As remarked in the introduction, these results have since
been strengthened considerably in \cite{BR07}, and indeed a strong law of large numbers was already implicit in \cite{Zern05}.  One of the key
purposes of this paper is to obtain analytically tractable
formulae for the coefficients $\theta(\beta,d)$ and $\Sigma(\beta,d)$
in the central limit theorem (see Section \ref{sec-discexp}),
allowing for a proof that $\beta\mapsto
\theta(\beta,d)$ is monotonically increasing. This result has since been
proved \ch{for $d\geq 9$} in \cite{HHmono08}, making crucial use of the methodology
in this paper.


\subsection{Random walk in a partially random environment}
\label{sec-RWRE}
In this section, we introduce a model of random walk in (partially) random environment (RWpRE).  The model we consider is a (nearest-neighbour, for simplicity) random walk in $\Z^d$, where $d=d_0+d_1$ with $d_1\ge 5$.  The random environment has the property that the random walker observed only when stepping in the coordinates $d_0+1,\dots, d$, behaves as a simple random walk in $d_1$ dimensions.   This is similar to the model studied in \cite{BSZ03}, but with two important differences.  Firstly, our results will only apply in the perturbative regime, where the transition probabilities are sufficiently close to their expected values.  Secondly we allow the probability of stepping in the $d_0+1,\dots, d$ coordinates to depend on the environment, provided that this probability is bounded away from 0, uniformly in the environment.  This situation is not allowed in \cite{BSZ03}, so our results for this model can be seen as a non-trivial extension to \cite{BSZ03} in the perturbative regime.

To be more precise about the model that we study, we require some additional notation.  Let $d=d_0+d_1\ge 6$ with $d_1\ge 5$.  Let $U_d$ be the set of unit vectors in $\Z^d$, and $\mc{P}(U_d)$ be the set of probability measures on $U_d$.  Let $W_{\cdot}=\{W_{\cdot}(u)\}_{u \in U_d}$ denote an element of $\mc{P}(U_d)$ and let $\mu$ be a probability measure on $\mc{P}(U_d)$,
satisfying the following:
\begin{itemize}
\item[(1)] The weight assigned to $U_{d_0}$ is not too large, i.e.~there exists some $\delta>0$ such that
\[\mu\Big(\sum_{u \in U_{d_0}}W_{\cdot}(u)\le 1-\delta\Big)=1.\]
\item[(2)] The weights assigned to $U_{d}\setminus U_{d_0}$ are ``fair", i.e.~for each $v \in U_d\setminus U_{d_0}$,
\[\mu\left(W_{\cdot}(v)=\frac{1-\sum_{u \in U_{d_0}}W_{\cdot}(u)}{2d_1}\right)=1.\]
\item[(3)] The weights cannot vary too much, i.e.~there exists some $\beta<1$ such that for each $u\in U_d$,
\eq
\lbeq{UE}\mu\left(|W_{\cdot}(u)-E_{\mu}[W_{\cdot}(u)]|<\beta\right)=1,
\en
where $E_{\mu}$ denotes expectation with respect to $\mu$.
\end{itemize}
Let $\nu$ be the product measure on $\mc{P}(U_d)^{\Z^d}$ obtained from $\mu$, i.e.~under $\nu$, $\{W_x\}_{x \in \Z^d}$ are independent with distribution $\ch{\mu}$.  The RWpRE in environment $W=\{W_x\}_{x \in \Z^d}$ is the Markov chain $\{X_n\}_{n\ge 0}$ such that $P_W(X_0=0)=1$, and $P_W(X_{n+1}=X_n+u|X_0,\dots,X_n)=W_{X_n}(u)$.  The annealed RWpRE is the (non-Markovian) random walk with law $\mQ$ obtained by averaging over all environments, i.e.
\[\mathbb{Q}(\vec{\omega}_n=\vec{x}_n)=\int P_W(\vec{X}_n=\vec{x}_n)d\nu.\]
The annealed transition probabilities are given by
    \eq
    \lbeq{prob_RWRE}
    p^{\vec{x}_i}(x_{i},x_{i+1})=\mQ(\omega_{i+1}=x_{i+1}|\vec{\omega}_i=\vec{x}_i)=\mE[W_{x_i}(x_{i+1}-x_i)|\vec{\omega}_i=\vec{x}_i].
    \en
Our main result for RWpRE is the following theorem, in which $\mQ_{\beta}$ denotes the above annealed law for a fixed (sufficiently small) choice of $\beta$.

\begin{theorem}[A CLT for RWpRE for $\ch{d_1>7}$]
\label{thm-1RWpRE}
Fix $\ch{d_1>7}$ and $d_0\ge 1$. Then, for every $\delta>0$ there exists $\beta_0 = \beta_0(d_1,d_0,\delta)>0$,
$\theta(\beta,\delta,d_1,d_0)$ and finite
$\Sigma=\Sigma(\beta,\delta,d_1,d_0)$ such that,
for all $\beta \leq \beta_0$,\\
(a)
\eq
\lbeq{meanRWpRE}
    \Ebold_\beta[\omega_n]
    =\theta n \big[1+O\big(\frac{1}{n}\big)\big].
\en
(b)
\eq
\lbeq{varRWpRE}
    \Var_\beta(\omega_n)
    = \Sigma n \big[1+O\big(\frac{\log n}{n^{1 \wedge \frac{d_1-6}{2}}}\big)\big].
\en
(c) $\omega_n$ satisfies a central limit theorem under
$\Qbold_\beta$, that is,
\eq
\lbeq{CLTRWpRE}
    \frac{\omega_n-\theta n}{\sqrt{n}} \convd {\cal N}(0,\Sigma).
\en
%
\end{theorem}


\begin{theorem}[A LLN for RWpRE for $\ch{d_1>4}$]
\label{thm-2RWpRE}
Fix $\ch{d_1>4}$ and $d_0\ge 1$. Then for each $\delta$, there exists $\beta_0 = \beta_0(d_1,d_0,\delta)>0$ and
$\theta=\theta(d_1,d_0,\delta)$ such that,
for all $\beta \leq \beta_0$,\\
(a)
\eq
\lbeq{meanRWpRELLN}
    \Ebold_\beta[\omega_n]
    =\theta n [1+O(\frac{\log{n}}{n^{1 \wedge \frac{d_1-4}{2}}})].
\en
(b) $\omega_n$ satisfies a law of large numbers under
$\Qbold_\beta$, that is,
\eq
\lbeq{LLNRWpRE}
    \frac{\omega_n}{n} \convp{\Qbold_\beta} \theta.
\en
%
\end{theorem}
In the above theorems, the $\beta_0$ arising from our analysis can be taken larger as $\delta$ increases.

Although results of a similar nature appear in \cite{BSZ03} and some of the references therein, we believe that this is a new result.  In particular, we do not assume that the random components of the environment are isotropic, nor that they have mean zero, nor that the random walker is transient in any particular direction.  However, as in \cite{BSZ03}, our analysis relies heavily on the fact that simple random walk in $d_1$ dimensions is, \ch{loosely} speaking, {\em very transient}.

\subsection{Overview of the method}
\label{sec-overview}
The main tool used is a {\it perturbation expansion for the two-point function}.
Such an expansion is often called a {\it lace expansion},
and takes the form of a recurrence relation
    \eq
    \lbeq{laceexp}
    c_{n+1}(x)= \sum_y D(y) c_n(x-y) + \sum_y \sum_{m=2}^{n+1} \pi_m(y) c_{n+1-m}(x-y)
    \en
for certain expansion coefficients $\{\pi_m\}_{m=2}^{\infty}$, and where
    \eq
    \lbeq{Ddef}
\ch{    D(x)=p^{\sss o}(o,x)}
    \en
is the transition probability function for the first step. A recurrence relation
such as \refeq{laceexp} is derived for the oriented percolation and
self-avoiding walk two-point functions, and plays an essential part in the proofs that these models are Gaussian
above the upper-critical dimension.  For self-avoiding walk, $c_n(x)$ equals the number of
$n$-step self-avoiding walks starting at 0 and ending at $x$, and $\sum_x c_{n}(x)$ equals the total number of self-avoiding walks, which grows
exponentially at a certain rate that needs to be determined in
the course of the proof.  For self-interacting random walks, $\sum_x c_{n}(x)=1$.  This essential difference gives rise
to a difference in the strategy for proofs.

In any lace expansion analysis, there are three main steps. The first is
the expansion in \refeq{laceexp}, which, for general self-interacting random
walks, will be derived in Section \ref{sec-lace}. The second step is to derive
bounds on the lace expansion coefficients.  These bounds will be derived in
Section \ref{sec-pibds}.  The final step is the analysis of the recurrence relation,
using the bounds on the lace expansion coefficients. For this analysis, we will make use
of induction.  The inductive analysis in this paper is intended for the perturbative regime (sufficiently small $\beta$), and is similar to the one in \cite{Hofs01}, where
a lace expansion was used to prove ballistic behaviour and a central limit theorem for
general one-dimensional weakly self-avoiding walk models.  In turn, this induction
was inspired by the analyses in \cite{HofHolSla98, HofSla02}.

In the induction argument,
we shall make use of the characteristic function of the end-point
of the $n$-step self-interacting random walk, which is the Fourier transform
    \eq
    \lbeq{cnk}
    \hat{c}_n(k) =\sum_{x\in \Z^d} e^{ik\cdot x} c_n(x).
    \en
Taking the Fourier transform of \refeq{laceexp} yields
    \eq
    \lbeq{fourierexp}
    \hat c_{n+1}(k)=\hat D(k)\hat c_{n}(k)+\sum_{m=2}^{n+1}\hat \pi_m(k)\hat c_{n+1-m}(k)
    \en

We shall present two separate induction arguments. The first
proves a law of large numbers as in Theorem \ref{thm-2ERW}
under relatively weak assumptions on the expansion coefficients, the
second is a more involved induction argument proving the
central limit theorem as in Theorems \ref{thm-1OERRW} and \ref{thm-1ERW}
under stronger assumptions on the expansion coefficients.

The remainder of the paper is organised as follows.
In Section \ref{sec-lace}, we present the expansion for
self-interacting random walks, which applies in the general
context described in Section \ref{sec-results}.  We also establish the formulae for the limiting speed and variance of the endpoint of the walk, assuming that these quantities exist. In
Sections \ref{sec-LLN} and \ref{sec-CLT}, we describe the induction arguments for the law of large numbers and central limit theorem respectively.  In Section \ref{sec-pibds}, we prove the bounds
on the lace expansion coefficients for the two models
under consideration. {\extend{In Section \ref{sec-expvar}, we prove the formula
for the variance stated in Theorem \ref{thm-var-form} in Section \ref{sec-lace}.
}}


\subsection{Recent applications of this method}
\label{sec-applic}
In this paper we have concentrated on deriving the expansion and on obtaining laws of large numbers and central limit theorems under strong conditions on the expansion coefficients.  However in Sections \ref{sec-expspeed} and \ref{sec-expvar} we also obtain formulae for the speed and variance, when these quantities are known to exist, under much weaker conditions on the expansion coefficients.

The speed of excited random walk is known to exist in all dimensions, e.g.~see \cite{BR07} and the results of this paper give a formula for that speed.  This formula is shown in \cite{HHmono08} to be monotone increasing in the excitation parameter in dimensions $d\ge 9$.  An excited random walk with opposing drift in a site-percolation cookie environment is studied in \cite{H08comp}.  A result of \cite{BSZ03} using cut-times and ergodicity shows that the speed of this model exists in high dimensions.  A formula for the annealed speed is then given by the results of this paper and it is shown in \cite{H08comp} that this formula is continuous in the excitation, percolation, and drift parameters, and strictly increasing in the excitation and percolation parameters in high dimensions.  In high dimensions, for each value of the drift parameter one can then establish phase transitions in the speed as one increases the percolation and excitement parameters.  In \cite{H08rwre}, certain models of random walk in i.i.d.~random environment, where at each site either the left or right step is not available, are studied in high dimensions.  In these models the existence of the speed is given by \cite{BSZ03}, a formula is provided by this paper, and it is possible to prove monotonicity of the speed as a function of the probability $p$ that the right step is available at the origin.

\section{The expansion for self-interacting random walks}
\label{sec-lace}
In this section, we perform and discuss the expansion for interacting random walks.
In Section \ref{sec-derexp}, we derive the expansion in \refeq{laceexp},
in Section \ref{sec-discexp}, we discuss the consequences of our expansion,
and in \ch{Sections \ref{sec-expspeed} and \ref{sec-expvar}, respectively,}
we identify the speed and variance from
our expansion formula, assuming that they exist and that the expansion formulae converge.

\subsection{Derivation of the expansion}
\label{sec-derexp}
Before we can start to prove \refeq{laceexp}, we need some more notation.
We will make use of the convolution of functions, which is
defined for absolutely summable
functions $f,g$ on $\Zd$ by
    \eq
    (f * g)(x) = \sum_y f(y) g(x-y),
    \en
so that we can rewrite \refeq{laceexp} as
    \eq
    \lbeq{laceexp2}
    c_{n+1}(x)= (D*c_n)(x) +\sum_{m=2}^{n+1} (\pi_m*c_{n+1-m})(x).
    \en

\ch{If $\vec{\eta}$ and $\vec{x}$ are two paths of length at least $j$ and $m$ respectively and such that $\eta_j=x_0$, then the concatenation $\ewalk_j\circ \vec{x}_m$ is defined by
 \eq
 \lbeq{concat}
    (\ewalk_j\circ \vec{x}_m)_{i}=\left\{
    \begin{array}{lll}
    &\eta_i&\text{when }0\le i\leq j,\\
    &x_{i-j} &\text{when }j \leq i \leq m+j.
    \end{array}\right.
    \en
Given $\ewalk_m$, we define a probability measure $\Qbold^{\ewalk_m}$ on walk paths starting from $\eta_m$, by specifying its value on particular cylinder sets (in a consistent manner) as follows
    \eq
    \Qbold^{\ewalk_m} (\walk_n=(x_0, x_1, \ldots, x_n))
    \equiv\prod_{i=0}^{n-1} p^{\ewalk_m\circ\vec{x}_i}(x_{i},x_{i+1}),
    \en
and extending the measure to all finite-dimensional cylinder sets in the natural (consistent) way.
We write $\Ebold^{\ewalk_m}$ for the expected value with respect to $\Qbold^{\ewalk_m}$, and define
    \eq
    \lbeq{cnhistory}
    c_{n}^{\ewalk_m}(\eta_m, x)= \Qbold^{\ewalk_m}
    (\omega_n=x).
    \en
}
Any path of length $n+1$ is a path of length $1$ concatenated with a path of
length $n$, so that, in terms of the above notation, we can
use \refeq{Ddef} to rewrite
    \eq
    \lbeq{exp1}
    c_{n+1}(x) =
    \sum_{\walkvec{0}{1}} D(\walkcoor{0}{1})\shift
    \sum_{\walkvec{1}{n}: \walkcoor{0}{1}\rightarrow x}\prod_{i=0}^{n-1}
    p^{\walkvec{0}{1}\circ \walkvec{1}{i}}(\walkcoor{1}{i},\walkcoor{1}{i+1}).
    \en
If we had $p^{\walkvec{0}{1}\circ \walkvec{1}{i}}=p^{\walkvec{1}{i}}$ \ch{ for all $\walkvec{0}{1}\circ \walkvec{1}{i}$},
then we would be back in the random walk case, since we would arrive at
    \eq
    \lbeq{RW}
    c_{n+1}(x) =\sum_{\walkvec{0}{1}} D(\walkcoor{0}{1})\shift
    \sum_{\walkvec{1}{n}: \walkcoor{0}{1}\rightarrow x}  \prod_{i=0}^{n-1}
    p^{\walkvec{1}{i}}(\walkcoor{1}{i},\walkcoor{1}{i+1})=(D*c_{n})(x).
    \en
For interacting random walks, $p^{\walkvec{0}{1}\circ \walkvec{1}{i}}$ does
not equal $p^{\walkvec{1}{i}}$ in general, and we are left to deal with the difference between the two.
For given \ch{$\ewalk_m$ and $\vec{x}_i$ we can }write
    \eq
    \lbeq{pmdiff}
\ch{    p^{\ewalk_m\circ \vec{x}_i}(x_i,x_{i+1})
    =
    p^{\vec{x}_i}(x_i,x_{i+1})+\big(p^{\ewalk_m\circ
    \vec{x}_i}-p^{\vec{x}_i}\big)(x_i,x_{i+1}).}
    \en
With this substitution, we have that
    \eq
    \lbeq{prodsub}
    \prod_{i=0}^{n-1}
    p^{\walkvec{0}{1}\circ \walkvec{1}{i}}(\walkcoor{1}{i},\walkcoor{1}{i+1})
    =\prod_{i=0}^{n-1}
    \big[p^{\walkvec{1}{i}}(\walkcoor{1}{i},\walkcoor{1}{i+1})+
    \big(p^{\walkvec{0}{1}\circ \walkvec{1}{i}}(\walkcoor{1}{i},\walkcoor{1}{i+1})-
    p^{\walkvec{1}{i}}(\walkcoor{1}{i},\walkcoor{1}{i+1})\big)\big].
    \en
In \refeq{prodsub}, the first term has `forgotten' the first step, while the
second term makes up for this mistake.  We would like to expand out the product in
\refeq{prodsub}.  Note that for all $\{a_i\}_{i=0}^{n-1}$ and $\{b_i\}_{i=0}^{n-1}$,
    \eq
    \lbeq{incl/excl}
    \prod_{i=0}^{n-1} (a_i+b_i)
    =\prod_{i=0}^{n-1} a_i+\sum_{j=0}^{n-1}
    \big(\prod_{i=0}^{j-1} (a_i+b_i) \big) b_j \big(\prod_{i=j+1}^{n-1} a_i\big),
    \en
where the empty products arising in $\prod_{i=0}^{j-1} (a_i+b_i)$ when $j=0$ and
$\prod_{i=j+1}^{n-1} a_i$ when $j=n-1$, are defined to be equal to 1.
Applying this to \refeq{exp1} with
    \[
    \ch{a_i=p^{\walkvec{1}{i}}(\walkcoor{1}{i},\walkcoor{1}{i+1}),
    \qquad
    b_i=\big(p^{\walkvec{0}{1}\circ \walkvec{1}{i}}-p^{\walkvec{1}{i}}\big)(\walkcoor{1}{i},\walkcoor{1}{i+1}),
    }
     \]
we arrive at
    \eqalign
    \lbeq{exp2}
    c_{n+1}(x) &=
    \sum_{\walkvec{0}{1}} D(\walkcoor{0}{1})\shift
    \sum_{\walkvec{1}{n}: \walkcoor{0}{1}\rightarrow x}\prod_{i=0}^{n-1}
    p^{\walkvec{1}{i}}(\walkcoor{1}{i},\walkcoor{1}{i+1})\nn\\
    &\qquad +\sum_{j=0}^{n-1}
    \sum_{\walkvec{0}{1}} D(\walkcoor{0}{1})\shift
    \sum_{\walkvec{1}{n}: \walkcoor{0}{1}\rightarrow x}
    \Big[\prod_{i=0}^{j-1}p^{\walkvec{0}{1}\circ \walkvec{1}{i}}(\walkcoor{1}{i},\walkcoor{1}{i+1})\Big]
    \big(p^{\walkvec{0}{1}\circ \walkvec{1}{j}}-p^{\walkvec{1}{j}}\big)(\walkcoor{1}{j},\walkcoor{1}{j+1})\nn\\
    &\qquad \times\Big[\prod_{i=j+1}^{n-1}
    p^{\walkvec{1}{i}}(\walkcoor{1}{i},\walkcoor{1}{i+1})\Big].
    \enalign
The first term equals $(D*c_{n})(x)$ by \refeq{RW}.
To rewrite the second term, we need some more notation.
We abbreviate
    \eq
    \Delta^{\smallsup{1}}_{j+1}= \big(p^{\walkvec{0}{1}\circ
    \walkvec{1}{j}}-p^{\walkvec{1}{j}}\big)(\walkcoor{1}{j},\walkcoor{1}{j+1}),
    \en
so that \refeq{exp2} becomes
    \eqalign
    c_{n+1}(x) &=
    (D*c_{n})(x)+\sum_{j=0}^{n-1}
    \sum_{\walkvec{0}{1}} D(\walkcoor{0}{1})\shift\shift
    \sum_{\walkvec{1}{n}: \walkcoor{0}{1}\rightarrow x}
    \Big[\prod_{i=0}^{j-1}p^{\walkvec{0}{1}\circ \walkvec{1}{i}}(\walkcoor{1}{i},\walkcoor{1}{i+1})\Big]
    \Delta^{\smallsup{1}}_{j+1}\Big[\prod_{i=j+1}^{n-1}
    p^{\walkvec{1}{i}}(\walkcoor{1}{i},\walkcoor{1}{i+1})\Big]\nn\\
    &= (D*c_{n})(x)+\sum_{j=0}^{n-1}
    \sum_{\walkvec{0}{1}} D(\walkcoor{0}{1})\sum_{\walkvec{1}{j+1}:\walkcoor{1}{0}=\walkcoor{0}{1}}\Big[\prod_{i=0}^{j-1}p^{\walkvec{0}{1}\circ \walkvec{1}{i}}(\walkcoor{1}{i},\walkcoor{1}{i+1})\Big]
    \Delta^{\smallsup{1}}_{j+1}\nn\\
    &\quad \qquad\qquad \times \shift\shift\sum_{(\walkcoor{1}{j+2}, \ldots, \walkcoor{1}{n}): \walkcoor{1}{n}=x}\Big[\prod_{i=j+1}^{n-1}p^{\walkvec{1}{i}}(\walkcoor{1}{i},\walkcoor{1}{i+1})\Big]
    .\lbeq{exp3}
    \enalign
From \refeq{cnhistory}, we have that
    \eq
    \sum_{(\walkcoor{1}{j+2}, \ldots, \walkcoor{1}{n}): \walkcoor{1}{n}=x}
    \Big[\prod_{i=j+1}^{n-1}
    p^{\walkvec{1}{i}}(\walkcoor{1}{i},\walkcoor{1}{i+1})\Big]
    =c_{n-j-1}^{\walkvec{1}{j+1}}(\walkcoor{1}{j+1},x).
    \en
Therefore, \refeq{exp3} is equal to
    \eqalign
    \lbeq{exp4}
    c_{n+1}(x)&=  (D*c_{n})(x)+\sum_{j=0}^{n-1}\sum_{\walkvec{0}{1}} D(\walkcoor{0}{1})\sum_{\walkvec{1}{j+1}}\Qbold^{\walkvec{0}{1}}(\vec{\omega}_j=\walkvec{1}{j})\Delta^{\smallsup{1}}_{j+1}~
    c_{n-j-1}^{\walkvec{1}{j+1}}(\walkcoor{1}{j+1},x)
    \enalign

For the second step of the expansion, we note that a type of two-point function $c_{n-j-1}^{\walkvec{1}{j+1}}(\walkcoor{1}{j+1},x)$ appears on the right side of \refeq{exp4}.  The second step of the expansion involves expanding out the dependence of this two-point function on the history $\walkvec{1}{j+1}$.  Given $\walkvec{1}{j+1}$ we write
\eqalign
\lbeq{2exp1}
c_{n-j-1}^{\walkvec{1}{j+1}}(\walkcoor{1}{j+1},x)=c_{n-j-1}(\walkcoor{1}{j+1}, x)
+\left(c_{n-j-1}^{\walkvec{1}{j+1}}(\walkcoor{1}{j+1},x)-c_{n-j-1}(\walkcoor{1}{j+1},x)\right).
\enalign
The contribution to \refeq{exp4} from the first term on the right of \refeq{2exp1} is
\eqalign
\lbeq{2exp1b}
\sum_{j=0}^{n-1}\sum_y\left[\sum_{\walkvec{0}{1}} D(\walkcoor{0}{1})\sum_{\walkvec{1}{j+1}}\Qbold^{\walkvec{0}{1}}(\vec{\omega}_j=\walkvec{1}{j})\Delta^{\smallsup{1}}_{j+1}I_{\{\omega^{(1)}_{j+1}=y\}}\right]~
    c_{n-j-1}(x-y)\equiv \sum_{m=2}^{n+1}\sum_y\pi_m^{\smallsup{1}}(y)c_{n+1-m}(x-y),
\enalign
where, for $m\ge 2$,
 \eqalign
    \lbeq{pi1def}
    \pi_m^{\smallsup{1}}(y)&=\sum_{\walkvec{0}{1}} D(\walkcoor{0}{1})\sum_{\walkvec{1}{m-1}}\Qbold^{\walkvec{0}{1}}(\vec{\omega}_{m-2}=\walkvec{1}{m-2})\Delta^{\smallsup{1}}_{m-1}I_{\{\omega^{(1)}_{m-1}=y\}}
    \enalign
To investigate the contribution to \refeq{exp4} from the term in brackets on the right of \refeq{2exp1}, we consider  the difference between $c_{n}^{\ewalk_m}(\eta_m, x)$ and $c_{n}(\eta_m, x)$
for general $\ewalk_m$, $n$ and $x$.
We first write
    \eq
    c_{n}^{\ewalk_m}(\eta_m,x)=\sum_{\walk_n^*:\eta_m\rightarrow x}
    \prod_{i=0}^{n-1} p^{\ewalk_m\circ\walk_i^*}(\omega_{i}^*,\omega^*_{i+1}),
    \en
and then use \refeq{pmdiff} and \refeq{incl/excl} to end up with
    \eq
    c_{n}^{\ewalk_m}(\eta_m, x)=c_{n}(\eta_m, x)
    +\sum_{j=0}^{n-1}\sum_{\walk_n^*:\eta_m\rightarrow x}
    \Big[\prod_{i=0}^{j-1} p^{\ewalk_m\circ\walk_i^*}(\omega^*_{i},\omega^*_{i+1})\Big]
    \big(p^{\ewalk_m\circ
    \walk^*_{j}}-p^{\walk^*_{j}}\big)(\omega^*_{j},\omega^*_{j+1})
    \prod_{i=j+1}^{n-1} p^{\walk^*_i}(\omega^*_{i},\omega^*_{i+1})\Big].
    \en
Therefore, similarly to \refeq{exp3}--\refeq{exp4}, we obtain
    \eq
    \lbeq{laceexpceta}
    c_{n}^{\ewalk_m}(\eta_m, x)=c_{n}(\eta_m, x)
    +\sum_{j=0}^{n-1}\sum_{\vec{\omega}_{j+1}^*}
    \Qbold^{\ewalk_m}(\vec{\omega}_{j}=\vec{\omega}^*_j)\Delta^{*}_{j+1}c_{n-j-1}^{\walk_{j+1}^*}(\omega_{j+1}^*,x)
    .
    \en
In \refeq{laceexpceta}, the first term is a regular two-point function, i.e., it does not depend on the history $\ewalk_m$.  In the correction term a history-dependent two-point function $c_{n-j-1}^{\walk_{j+1}^*}$ appears to which we can iteratively use \refeq{laceexpceta}.
Thus, with $m=j+2$,
    \eqalign
    \lbeq{exp5}
    c_{n+1}(x)=&(D*c_{n})(x)+\sum_{m=2}^{n+1}(\pi_m^{\smallsup{1}}*c_{n-m+1})(x)\\
    &+\sum_{j_1, j_2}\sum_{\walkvec{0}{1}}D(\walkcoor{0}{1})\sum_{\walkvec{1}{j_1+1}}\sum_{\walkvec{2}{j_2+1}}
    \Qbold^{\walkvec{0}{1}}(\vec{\omega}_{j_1}=\walkvec{1}{j_1})\Delta^{\smallsup{1}}_{j_1+1}\Qbold^{\walkvec{1}{j_1+1}}(\vec{\omega}_{j_2}=\walkvec{2}{j_2})\Delta^{\smallsup{2}}_{j_2+1}
    c_{n-j_1-j_2-2}^{\walkvec{2}{j_2+1}}(\omega^{\sss (2)}_{j_2+1},x),\nn
    \enalign
where we write, for $N\geq 1$,
    \eq
    \lbeq{DeltaNdef}
    \Delta^{\smallsup{N}}_{j_{N}+1}= \big(p^{\walkvec{N-1}{j_{N-1}+1}\circ
    \walkvec{N}{j_{N}}}-p^{\walkvec{N}{j_{N}}}\big)(\walkcoor{N}{j_{N}},\walkcoor{N}{j_{N}+1}),
    \en
with $j_0\equiv 0$.

For $N\geq 1$, we let $\mc{A}_{m,\sN}=\{\vec{j}\in \Z_+^N:j_1+ \cdots+ j_{\sss N}=m-N-1\}$ and further define
    \eqalign
    \pi_m^{\smallsup{N}}(y)=&\sum_{\vec{j}\in \mc{A}_{m,N}}\sum_{\walkvec{0}{1}}\sum_{\walkvec{1}{j_{\sss 1}+1}}\dots\sum_{\walkvec{N}{j_{\sss N}+1} }I_{\{\omega^{(N)}_{j_{\sN}+1}=y\}}D(\walkcoor{0}{1})
    \prod_{n=1}^{N}\Delta^{\sss (n)}_{\sss j_n+1}\prod_{i_{n}=0}^{j_{n}-1}p^{\walkvec{n-1}{j_{n-1}+1}\circ \walkvec{n}{i_{n}}}\left(\walkcoor{n}{i_{n}},\walkcoor{n}{i_{n}+1}\right)   \lbeq{piNdef}\\
 =&   \sum_{\vec{j}\in \mc{A}_{m,N}}\sum_{\walkvec{0}{1}}\sum_{\walkvec{1}{j_{\sss 1}+1}}\dots\sum_{\walkvec{N}{j_{\sss N}+1} }I_{\{\omega^{(N)}_{j_{\sN}+1}=y\}}D(\walkcoor{0}{1})
    \prod_{n=1}^{N}\Delta^{\sss (n)}_{\sss j_n+1}\mQ^{\walkvec{n-1}{j_{n-1}+1}}(\vec{\omega}_{j_n}=\walkvec{n}{j_n}).\nn
   \enalign
which is zero when $N+1>m$.  Note that \refeq{piNdef} reduces to \refeq{pi1def} in the case $N=1$.  Then define
    \eq
    \pi_m(y)=\sum_{N=1}^{\infty} \pi_m^{\smallsup{N}}(y).
    \en
 We emphasize that, conditionally on $\walkvec{M}{j_{M}+1}$,
the probability measure ${\mathbb Q}_{\sss M+1}^{\walkvec{M}{j_{M}+1}}$ is the law of
$\walkvec{M+1}{j_{M+1}+1}$, i.e., that $\walkvec{M}{j_{M}+1}$ acts as the history for $\walkvec{M+1}{j_{M+1}+1}$.

Equation \refeq{laceexp} follows by iteratively replacing the two-point function in \refeq{laceexpceta} by using the equality \refeq{laceexpceta}, until the second term on the right of \refeq{laceexpceta} vanishes.   This must happen when $N=n+1$.
This completes the derivation of the expansion.

\subsection{Discussion of the expansion}
\label{sec-discexp}
In this section, we discuss the consequences of the expansion in \refeq{laceexp}.

\paragraph{The lace expansion coefficients.}
The lace expansion coefficients involve the factors
    \eq
    \Delta^{\smallsup{N}}_{j_{N}+1}= \big(p^{\walkvec{N-1}{j_{N-1}+1}\circ
    \walkvec{N}{j_{N}}}-p^{\walkvec{N}{j_{N}}}\big)(\walkcoor{N}{j_{N}},\walkcoor{N}{j_{N}+1})
    \en
in \refeq{DeltaNdef}. This difference is identically zero when the histories
$\walkvec{N-1}{j_{N-1}+1}\circ \walkvec{N}{j_N}$ and $\walkvec{N}{j_N}$ give the same transition
probabilities to go from $\walkcoor{N}{j_N}$ to $\walkcoor{N}{j_N+1}$.  For excited random walk, $\Delta^{\smallsup{N}}_{j_N+1}$ is non-zero precisely when $\walkcoor{N}{j_N}$ has already been visited by $\walkvec{N-1}{j_{N-1}+1}$, but not by $\walkvec{N}{j_{N}-1}$, so that
\begin{equation}
\begin{split}
\label{e:DeltaExRW}
|\Delta^{\smallsup{N}}_{j_{N}+1}|
    \leq &|\Delta^{\smallsup{N}}_{j_{N}+1}|I_{\{\walkcoor{N}{j_{N}}\in \walkvec{N-1}{j_{N-1}+1}\}}I_{\{\walkcoor{N}{j_{N}}\notin \walkvec{N}{j_{N}-1}\}}\\
    \leq &C\beta
    I_{\{\walkcoor{N}{j_{N}}\in \walkvec{N-1}{j_{N-1}+1}\}}I_{\{\walkcoor{N}{j_{N}}\notin \walkvec{N}{j_{N}-1}\}}\le C\beta I_{\{\walkcoor{N}{j_{N}}\in \walkvec{N-1}{j_{N-1}}\}}.
\end{split}
\end{equation}

\begin{figure}
\includegraphics[scale=.75]{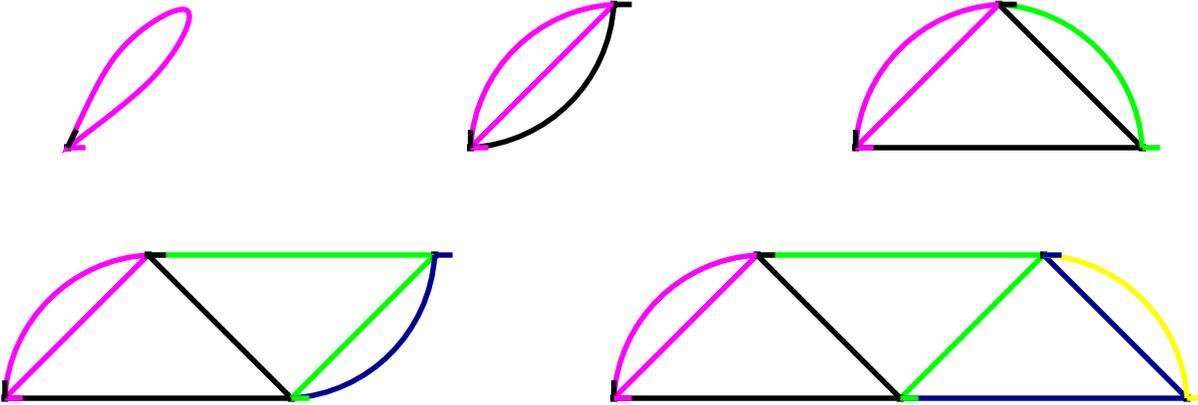}
\caption{The diagrams for $\pi_m^{\smallsup{N}}$, $N=1, \dots, 5$, arising from the expansion and the bound \refeq{DeltaERRW} for both models.  The subwalks (indicated by different shades) in the diagrams have the previous subwalk as their history.  An intersection of two subwalks and a small factor $\beta$ appears at each vertex.}
\label{fig:expdiags}
\end{figure}

For once-edge-reinforced
random walk, the difference \refeq{DeltaNdef} is nonzero exactly when the vertex $\walkcoor{N}{j_N}$ has already been visited by $\walkvec{N-1}{j_{N-1}+1}$ via an edge that was not traversed by $\walkvec{N}{j_N}$. Therefore,
we also have for once-edge-reinforced random walk that
    \eq
    \lbeq{DeltaERRW}
    |\Delta^{\smallsup{N}}_{j_{N}+1}|\leq C
    \beta
    I_{\{\walkcoor{N}{j_{N}}\in \walkvec{N-1}{j_{N-1}}\}}.
    \en
\ch{For RWpRE, a similar bound holds as follows.  From \refeq{prob_RWRE},
\eq
\Delta^{\smallsup{N}}_{j_{N}+1}=\mE\big[W_{\walkcoor{N}{j_{N}}}(\walkcoor{N}{j_{N}+1}-\walkcoor{N}{j_{N}})\big|\vec{\omega}_{j_{N-1}+1+j_{N}}=\walkvec{N-1}{j_{N-1}+1}\circ    \walkvec{N}{j_{N}}\big]-\mE\big[W_{\walkcoor{N}{j_{N}}}(\walkcoor{N}{j_{N}+1}-\walkcoor{N}{j_{N}})\big|\vec{\omega}_{j_{N}}=\walkvec{N}{j_{N}}\big].\nn
\en
By definition the random environment is site-wise independent, so the only information about $W_{\walkcoor{N}{j_{N}}}$ contained in the history of the path is in the departures from the site $\walkcoor{N}{j_{N}}$.  Trivially every departure from $\walkcoor{N}{j_{N}}$ by $\walkvec{N}{j_{N}}$ is also a departure from $\walkcoor{N}{j_{N}}$ by $\walkvec{N-1}{j_{N-1}+1}\circ \walkvec{N}{j_{N}}$, and any additional departures from this site by $\walkvec{N-1}{j_{N-1}+1}\circ \walkvec{N}{j_{N}}$ are actually departures from $\walkcoor{N}{j_{N}}$ by $\walkvec{N-1}{j_{N-1}}$.  Thus $\Delta^{\smallsup{N}}_{j_{N}+1}$ is non-zero only if $\walkcoor{N}{j_{N}}\in \walkvec{N-1}{j_{N-1}}$.  It then follows immediately from \refeq{UE} that
\eq
\lbeq{DeltaRWRE}
|\Delta^{\smallsup{N}}_{j_{N}+1}|\leq 2
    \beta
    I_{\{\walkcoor{N}{j_{N}}\in \walkvec{N-1}{j_{N-1}}\}}.
    \en}
We conclude that for \ch{all models under consideration}, each factor $|\Delta^{\smallsup{N}}_{j_{N}+1}|$:
\begin{enumerate}
    \item enforces an intersection between the path and its previous history;
    \item gives rise to a factor $\beta$, making $\pi_m^{\smallsup{N}}(y)$ small when $\beta$ is sufficiently small and $N$ is large.
\end{enumerate}
The quantities $\pi_m^{\smallsup{N}}(y)$ combined with the bound \refeq{DeltaERRW} for both models, can be represented by diagrams of the form displayed in Figure \ref{fig:expdiags} for $N=1, \dots, 5$.  The first step is special, as it has no history.  Thereafter, each subwalk $\walkvec{i}{j_i+1}$ (indicated by shading in Figure \ref{fig:expdiags}) has the previous subwalk $\walkvec{i-1}{j_{i-1}+1}$ as its history.  The apparent similarity with the self-avoiding walk diagrams (see for example \cite{HarSla92b}) is natural due to the intersections enforced by the factors $\Delta^{\smallsup{i}}_{j_{i}+1}$ as described above.  A small factor $\beta$ arises from each intersection (represented by vertices in Figure \ref{fig:expdiags}), and the number of intersections increases with the complexity of the diagram.

\paragraph{The speed and variance.}
By convention our vectors are considered to be column vectors.   Thus if $\theta\in \R^d$, then $\theta\theta^t$ is a $d\times d$ matrix with real entries.

The limiting speed $\theta=\theta(\beta,d)$ and
covariance matrix $\Sigma=\Sigma(\beta,d)$ appearing in Theorems
\ref{thm-1OERRW}--\ref{thm-2ERW} are given by
    \eqarray
    \theta(\beta,d)& = &\zerospeed - i\sum_{m=2}^\infty \nabla\hat{\pi}_m(0),\lbeq{theta}\\
    \Sigma(\beta,d)& = &\zerovar-\theta\theta^t-\sum_{m=2}^\infty
        \nabla^2\Big[e^{-i\theta\cdot k(m-1)}\hat{\pi}_m(k)\Big]_{k=0},\lbeq{sigma}
    \enarray
where $\zerospeed$ is the expected drift of the transition probability $D=\zerop$,
i.e.,
    \eq
    \lbeq{zerospeeddef}
    \zerospeed=\sum_{x\in \Z^d} x D(x),
    \en
while $\zerovar$ is the covariance matrix of $D=\zerop$ given by
    \eq
    \lbeq{zerovardef}
    (\zerovar)_{i,j}= \sum_{x\in \Z^d} x_ix_j D(x),
    \en
and $\nabla f(k)$ is the vector of derivatives of $k\mapsto f(k)$, while $\nabla^2 f(k)$ is the matrix consisting of the double derivatives of
$k\mapsto f(k)$.

These formulas can be heuristically derived from the recurrence relation
\refeq{laceexp}. Indeed, take the Fourier transform to obtain
    \eq
    \lbeq{laceexpk}
    \hat{c}_{n+1}(k)= \hat D(k) \hat c_n(k) + \sum_{m=2}^{n+1} \hat\pi_m(k) \hat c_{n+1-m}(k).
    \en
Now replace $\hat{c}_{l}(k)$ throughout the recurrence relation by
$e^{i\theta\cdot k l-\hlf k^t\Sigma k l}$,
in accordance with Theorem \ref{thm-1OERRW}(c)--\ref{thm-1ERW}(c). Then,
dividing by $e^{i\theta\cdot k n-\hlf k^t\Sigma k n}$, we obtain
    \eq
    \lbeq{laceexpkb}
    e^{i\theta\cdot k-\hlf k^t\Sigma k}
    \approx \hat D(k)+ \sum_{m=2}^{n+1}
    \hat\pi_m(k)e^{-i\theta\cdot k (m-1)+\hlf k^t\Sigma k(m-1)}.
    \en
Expanding to linear order in $k$ yields \refeq{theta} and expanding to second order in $k$ yields \refeq{sigma}, when we note that $\Sigma$ (as defined in \refeq{var} and \refeq{varERW}) must be symmetric, and
    \eq
    \lbeq{Dhattaylor}
     \hat\pi_m(0)=0 \quad \text{and} \quad \hat{D}(k)=1+ik\cdot \zerospeed-\frac{1}{2}k^t\zerovar k +O(|k|^3).
    \en

The results in this paper, as well as the proofs, follow part of the ideas in
\cite{Hofs01}, where it was shown that certain weakly self-avoiding walk models
in $d=1$ behave ballistically.

\subsection{The formula for the speed}
\label{sec-expspeed}

In this section, we show that, when the speed is
proved elsewhere to exist, {\it and} our formula for the speed in
\refeq{theta} converges, then in fact
\refeq{theta} identifies the speed. For example, for ERW in dimensions $d=2,\ldots, 5$,
where Theorem \ref{thm-2ERW} does not apply, it is known (e.g. \cite{BR07}) that the speed exists almost surely.

\begin{theorem}[The speed formula]
\label{prp:speedmatch}
If $\lim_{n \ra \infty}\sum_{m=2}^n \sum_x x\pi_m(x)$ exists
and $n^{-1}\omega_n\convp{{\mathbb Q}_\beta} \theta$, then
    \eqalign
    \theta(\beta,d)& = \sum_x x p^{\sss 0}(0,x)
    +\sum_{m=2}^\infty \sum_x x\pi_m(x)\lbeq{theta-rep}.
    \enalign
\end{theorem}

\proof Multiplying \refeq{laceexp} by $x=y+(x-y)$, summing, and using the facts that
$\sum_x c_n(x)=\sum_x p^{\sss 0}(0,x)=1$ and $\sum_x \pi_m(x)=0$,
we obtain
    \eqalign
    \sum_x xc_{n+1}(x)=&\sum_y yp^{\sss 0}(0,y)+\sum_x xc_n(x)
    +\sum_{m=2}^{n+1}\sum_y y\pi_m(y).\lbeq{samespeed}
    \enalign
Now $\sum_x x c_{n}(x)=\mE[\omega_n]$, so rearranging \refeq{samespeed} we obtain
    \eqalign
    \lbeq{meandrift1}
    \mE[\omega_{n+1}-\omega_n]=\zerospeed+\sum_{m=2}^{n+1}\sum_y y\pi_m(y).
    \enalign
The right hand side converges if and only if the left hand side does.
Thus, under the assumption that $\lim_{n \ra \infty}\sum_{m=2}^n \sum_x x\pi_m(x)\equiv\tilde\theta(\beta,d)-\zerospeed$ exists,
we obtain that
    \eqn{
    \lbeq{meandrift2}
    \lim_{n\rightarrow \infty} \mE[\omega_{n+1}-\omega_n]=\tilde\theta(\beta,d).
    }
In turn, \refeq{meandrift2} implies that
    \eqn{
    \lim_{n\rightarrow \infty} \mE[n^{-1}\omega_n]=\tilde\theta(\beta,d).
    }
When  $n^{-1}\omega_n\convp{{\mathbb Q}_\beta} \theta(\beta,d)$, by bounded convergence
and the fact that $|\omega_n|\le nL$ since the maximal step size
of our self-interacting random walks is $L$, we have that
    \eqn{
    \lim_{n\rightarrow \infty} \mE[n^{-1}\omega_n]=\theta(\beta,d),
    }
so that, as required, $\theta(\beta,d)=\tilde\theta(\beta,d)$.
\qed

\subsection{The formula for the variance}
\label{sec-expvar}
In this section we prove a result about the variance of the endpoint of the walk, similar to that obtained above for the speed.
Define $a_m^{[i]}:=\sum_yy^{[i]}\pi_m(y)$. Then, we have the following formula for the variance
of self-interacting random walks in terms of the lace expansion coefficients:

\begin{theorem}[The variance formula]
\label{thm-var-form}
Suppose that for each $i,j\in \{1,2,\dots,d\}$,
\eq
\lbeq{lim_exists}
\lim_{n\ra \infty}\frac{\mE[\omega^{[i]}_{n}\omega^{[j]}_{n}]-\mE[\omega^{[i]}_{n}]\mE[\omega^{[j]}_{n}]}{n}= \Sigma_{ij}, \quad \text{ and }\quad \sum_{m=2}^{\infty}\sum_yy^{[i]}y^{[j]}\pi_m(y)<\infty,
\en
and that either
\begin{itemize}
\item[{\em (i)}]   $\mE[\omega_n]=0$ for each $n$, or
\item[{\em (ii)}]   $n^{-1}\omega_n\convp{{\mathbb Q}_\beta}\theta$, and
$\sum_{m=2}^{\infty}(m-1)|a_m^{[i]}|<\infty$.
\end{itemize}
Then
\eqalign
\lbeq{lim_covariance}
\Sigma_{ij}=(\zerovar)_{ij}-\theta^{[i]}\theta^{[j]}-\sum_{m=2}^{\infty}\left[\theta^{[i]}(m-1)a_m^{[j]}+\theta^{[j]}(m-1)a_m^{[i]}-\sum_y y^{[i]}y^{[j]}\pi_m(y)\right].
\enalign
\end{theorem}

The proof of Theorem \ref{thm-var-form} is an adaptation of that
of the speed formula in Theorem \ref{prp:speedmatch} above, and is
deferred to {\extend{Section \ref{sec-thm-var-form}}} {\short{\cite{HofHol10b}.}}

\section{Induction for the weak law of large numbers}
\label{sec-LLN}
In this section we prove a law of large numbers from the recurrence relation \refeq{laceexp}, or, more precisely, its Fourier transform \refeq{laceexpk}, assuming certain bounds on the coefficients $\hat{\pi}_m(k)$.  The bounds roughly correspond to upper bounds on the accuracy of the Taylor approximation of $\hat{\pi}_m(k)$ up to {\em first} order.

We start by formulating a general assumption (which must be verified for a specific model), and prove the main result, Theorem \ref{thm-LLNgen}, under this assumption.

\begin{ASSLLN}
\label{ass:LLN}
There exists a sequence $\{b_m\}_{m \ge 1}$, independent of $\beta$ and with $b_1\ge 1$, and a constant $\veta=\veta(d)$ satisfying $\lim_{\beta\rightarrow 0}\veta(d)=0$
 such that
    \eq
    \lbeq{pibdscfirst}
    \hat{\pi}_m(0)=0, \qquad |\nabla
    \hat{\pi}_m(0)|
    \leq \veta b_m,
    \qquad
    |\nabla^2\hat{\pi}_m(0)|\leq
    \veta m b_m,
    \en
and uniformly in $k \in [-\pi,\pi]^d$,
    \eq
    \lbeq{piLLNbound}
    |\hat{\pi}_m(k)|\leq \veta|k|b_m,
    \qquad |\hat{\pi}_m(k)-k\cdot\nabla\hat{\pi}_m(0)|\le \veta|k|^2mb_m,
    \en
where
    \eq
    \lbeq{bmsummable}
    B\equiv \sum_{m=1}^\infty b_m<\infty \quad \text{and}\quad  B'\equiv \sup_{n}\frac{(\log{(n\vee 3)})^2}{n}\sum_{m=1}^{n} mb_m<\infty.
    \en
\end{ASSLLN}

\begin{theorem}[Weak law of large numbers]
\label{thm-LLNgen}
When Assumption (LLN) holds, there exist $\beta_0=\beta_0(d)>0$
and $\theta=\theta(\beta)$ such that for all $\beta\le \beta_0$,
    \eq
    \lbeq{mean-LLN}
    \Ebold_\beta[\omega_n]
    = \theta n \Big[1+O\Big(\frac{1}{n}\sum_{m=1}^\infty (n\wedge m) b_m\Big)\Big].
    \en
Furthermore, there exists $C>0$ such that for every $k\in \R^d$,
    \eq
    \lbeq{LLN}
\log\left(\Ebold_\beta[e^{ik\cdot \omega_n/n}]\right) =
    ik\cdot \theta+O\left(\frac{|k|}{n}e^{C|k|}\sum_{m=1}^\infty (n\wedge m) b_m\right)
    +O\left(\frac{|k|^2}{n} \sum_{m=1}^n m b_m\right),
\en
where the constant $\theta$ given by \refeq{theta} is model dependent.
\end{theorem}
\begin{REM}
Observe that $n^{-1}\sum_{m=1}^{\infty}(n \wedge m)b_m=o(1)$ and $n^{-1}\sum_{m=1}^{n}mb_m=o(1)$ when \refeq{bmsummable} holds.  Thus \refeq{LLN} implies that $\lim_{n \rightarrow \infty}\Ebold_\beta[e^{ik\cdot \omega_n/n}]=e^{ik\cdot \theta}$, which is equivalent to the statement of convergence in probability, $\omega_n/n\convp{{\mathbb Q}_\beta}\theta$.
\end{REM}
\vskip0.5cm

Note that since $D$ has finite range, there exists a constant $C_1\ge 1$ independent of $\beta$ such that
    \eq
    \lbeq{DTaylor}
    |\hat{D}(k)-1-ik\cdot \zerospeed|\leq C_1|k|^2,
    \en
and let $K_1=2C_1$, which is independent of $\beta$.

We will frequently use the following
lemma, whose proof follows easily by applying Taylor's Theorem at $t=0$ to the map from $\R\ra \C$ given by $t \mapsto e^{tx}$:
\begin{lemma}
\label{lem-Taylor}
For all $x\in \mathbb{C}, j\in \mathbb{N}$,
$$
    \left|e^{x}-\sum_{l=0}^j \frac{x^l}{l!}\right|  \leq
    \frac{|x|^{j+1}}{(j+1)!} e^{|{\rm Re}(x)|},
$$
where ${\rm Re}(x)$ is the real part of $x$.
\end{lemma}

Set $\theta_1=\zerospeed$, and, for $n\ge 2$, we define the
following approximation to $\theta$:
    \eq
    \theta_n = \zerospeed - i\sum_{m=2}^n \nabla\hat{\pi}_m(0).\lbeq{thetan}
    \en

Our induction hypothesis for the law of large numbers
in Theorem \ref{thm-LLNgen}
is that the following bound holds for
all $\beta\leq \beta_0$, some $\delta<1$ independent of $\beta$ and all $0\leq j\leq n$:

For $|k|\le \delta\log(n\vee 3)/n$ and
some $K\ge 1$ independent of $\beta$ we can write,
    \eq
\quad \hat c_j(k) = \exp{\Big[\sum\limits_{l=1}^j \big(ik \cdot
    \theta_l+e_l(k)\big)\Big]}\quad \mbox{where
    }\quad
    |e_j(k)|\leq K|k|^2\sum_{l=1}^j lb_l,
    \lbeq{IHLLN}
    \en
where the empty sum, arising when $j=0$, is defined to be 0, and where,
for $n=0$, the equation is valid for all $k\in [-\pi,\pi]^d$.

The initialisation of the induction (the $n=0$ case) holds trivially since $1=e^0$.  In Section \ref{sec-LLNadv} we will advance the induction hypothesis.
In Section \ref{sec-pfLLN} we will use it to prove Theorem \ref{thm-LLNgen}.

\subsection{The LLN induction advanced}
\label{sec-LLNadv}
We fix $n\geq 0$.
The induction step will be achieved as
soon as we are able to write
    \eq
    \frac{\hat{c}_{n+1}(k)}{\hat{c}_n(k)}=\exp{\big[ik\cdot \theta_{n+1}+e_{n+1}(k)\big]},
    \en
for $e_{n+1}(k)$ satisfying the required bound. For this,
we write
    \eq \lbeq{e'}
    \frac{\hat{c}_{n+1}(k)}{\hat{c}_{n}(k)}
    =
    1+ik\cdot\theta_{n+1}+e'_{n+1}(k)
    \en
and then set
    \eq
    \lbeq{errors}
    e_{n+1}(k) =
    \log{\big[1+ik\cdot\theta_{n+1}+e'_{n+1}(k)\big]}
    - ik\theta_{n+1}.
    \en
The following lemma is a trivial consequence of \refeq{zerospeeddef} and \refeq{thetan}:
\begin{lemma}
\label{lem-thetabds}
We have $|\zerospeed|\le L$ and when Assumption (LLN) holds we have
$|\theta_n|\le L+\vep_{\beta} B$ for every $n$.
\end{lemma}

Let
    \eq
    \lbeq{Bn-def}
    B_n=\sum_{m=1}^n mb_m.
    \en
We note that by the second bound in \refeq{bmsummable}, and uniformly in
$k$ such that $|k|\le \delta n^{-1}\log{(n\vee 3)}$, we have
    \eq
    \lbeq{Bnk-bds}
    B_{n+1}|k| \leq \delta B', \qquad nB_n |k|^2 \leq \delta^2 B'.
    \en
These bounds will be frequently used in what follows.

Choose $\beta_0>0$ so that $\veta\le 1$ for all
$\beta\le \beta_0$, and suppose that the required bound
\refeq{IHLLN} holds for $e'_{n+1}(k)$ with constant $K_1$.
By Lemma \ref{lem-thetabds}, $|k||\theta_{n+1}|+|e'_{n+1}(k)|\le 1/2$
for $|k|\le \delta\log(n\vee 3)/n$ when $\delta\leq(2(L+B+K_1B')\log{3})^{-1}$.
Therefore we may apply Taylor's Theorem $|\log(1+x)-x|\le 4|x|^2$ for
$|x|\le 1/2$, to \refeq{errors}.  This implies that when the
required bound holds for $e'_{n+1}(k)$ with constant $K_1$, it
also holds for $e_{n+1}(k)$ for some $K$ independent of $\beta$,
since the terms of order $k$ in \refeq{errors} cancel.  Specifically, if
$|e'_{n+1}(k)|\le K_1|k|^2B_{n+1}$, then, using also \refeq{Bnk-bds} and
$(x+y)^2\leq 2x^2+2y^2$,
    \eqalign
    |e_{n+1}(k)|\le &4(|k||\theta_{n+1}|+|e'_{n+1}(k)|)^2+|e'_{n+1}(k)|\nn\\
    \le & 8|k|^2(L+B)^2+8K_1^2B_{n+1}^2|k|^4+K_1B_{n+1}|k|^2\le |k|^2(8(L+B)^2+8B_{n+1}K_1^2\delta B'+B_{n+1}K_1)\nn\\
    \le &KB_{n+1}|k|^2,
    \enalign
for $K\ge 8(L+B)^2+8\delta B'K_1^2+K_1,$ which is independent of $\beta$.

The rest of this section will be devoted to the proof of the
following lemma:
\begin{lemma}
\label{lem-ej'}
There exists $\beta_0$ such that for all $\beta\le \beta_0$, if $e_j(k)$ satisfies the bound in \refeq{IHLLN}
for all $j\leq n$ and $|k|\le \delta (n+1)^{-1}\log{((n+1)\vee 3)}$ then for such $k$,
    \eqarray
    |e'_{n+1}(k)|  &\leq &K_1B_{n+1}|k|^2.\lbeq{e'bound}
    \enarray
\end{lemma}
\proof
Divide the recursion relation \refeq{fourierexp} by
$\hat{c}_n(k)$ and use the equality $\hat{\pi}_m(0)=0$ of \refeq{pibdscfirst} to obtain
    \eqarray
    \frac{\hat{c}_{n+1}(k)}{\hat{c}_{n}(k)} & = &\hat{D}(k)+\sum_{m=2}^{n+1}
    [\hat{\pi}_m(k)-\hat{\pi}_m(0)]\frac{\hat{c}_{n+1-m}(k)}{\hat{c}_{n}(k)}.
    \lbeq{recf_n(k)}
    \enarray
We can rewrite \refeq{recf_n(k)} as
    \eqarray
    \frac{\hat{c}_{n+1}(k)}{\hat{c}_{n}(k)}
        & = &1~+ ik\cdot \theta_{n+1}+e'_{n+1}(k),\nn
    \enarray
where
    \eqarray
    \lbeq{edashLLN}
    e'_{n+1}(k) & =&[\hat{D}(k)-1-ik\cdot \zerospeed]
    +\sum_{m=2}^{n+1} \Big[\hat{\pi}_m(k)-
    k\cdot \nabla\hat{\pi}_m(0)\Big]+\sum_{m=2}^{n+1}
    \hat{\pi}_m(k)\big[\frac{\hat{c}_{n+1-m}(k)}{\hat{c}_n(k)}-1\big].\nn
    \enarray
The first term is taken care of by \refeq{DTaylor}.  Furthermore, by \refeq{piLLNbound}, we have that
    \eq
    \sum_{m=1}^{n+1}|\hat{\pi}_m(k)-
    k\cdot \nabla\hat{\pi}_m(0)|\leq \veta|k|^2 \sum_{m=1}^{n+1}mb_m=\veta B_{n+1}|k|^2.
    \en
Finally, using Lemma \ref{lem-Taylor} and the induction
hypothesis \refeq{IHLLN} for $e_l(k)$ with $l \le n$,
which is allowed since $|k|\le \delta\log((n+1)\vee 3)/(n+1)$
implies that also $|k|\le \delta\log(n\vee 3)/n$,
    \eqalign
    \Big|\frac{\hat{c}_{n+1-m}(k)}{\hat{c}_n(k)}-1\Big|
    =&
    \Big|\exp{\Big[-\sum\limits_{l=n+1-m}^n \big(ik\cdot \theta_l+e_l(k)\big)\Big]}
    -1\Big|\le m(|k|(L+B)+KB_{n+1}|k|^2)e^{mKB_{n+1}|k|^2}\nn\\
    \le& m(|k|(L+B)+KB_{n+1}|k|^2)e^{KB'\delta^2},
    \enalign
by the second inequality in \refeq{Bnk-bds}.

Using the first bound in \refeq{piLLNbound}, it follows that
    \eqalign
    \lbeq{piconvratiobound}
    \Big|\sum_{m=2}^{n+1}
    \hat{\pi}_m(k)\big[\frac{\hat{c}_{n+1-m}(k)}{\hat{c}_n(k)}-1\big]\Big|
    \leq &\sum_{m=2}^{n+1}|k|\veta b_m e^{K B'\delta^2}m(|k|(L+B)+B_{n+1}K|k|^2)\nn\\
    \leq &\veta|k|^2B_{n+1}e^{KB'\delta^2}\Big((L+B)+B_{n+1}K|k|\Big)\nn\\
    \le &\veta |k|^2B_{n+1}e^{KB'\delta^2}\big((L+B)+\delta B'K\big),
    \enalign
where in the last inequality we used the first bound in \refeq{Bnk-bds}.
Summarising \refeq{edashLLN}-\refeq{piconvratiobound} we have
    \eq
    |e'_{n+1}(k)|\leq C_1|k|^2+\veta B_{n+1}|k|^2
    +\veta e^{KB'\delta^2}\big((L+B)+\delta B'K\big)B_{n+1}|k|^2.
    \en
Recall that $K_1=2C_1$ and $B_{n+1}\geq b_1\geq 1$.
We choose $\beta_0$ sufficiently small so that both $\veta\le 1$ for all $\beta\le \beta_0$ and
    \eq
    \veta
    \Big(1+e^{KB'\delta^2}\big((L+B)+\delta B'K\big)\Big)\le \frac{K_1}{2}.
    \en

Then, we conclude that
$C_1+\veta B_{n+1}\big(1+e^{KB'\delta^2}\big((L+B)+\delta B'K\big)\big)
\le K_1B_{n+1}$ and therefore \refeq{e'bound} holds as required for
all $\beta\le \beta_0$.  This completes the proof of Lemma \ref{lem-ej'}.
\qed

\subsection{Proof of Theorem \ref{thm-LLNgen}}
\label{sec-pfLLN}
To prove \refeq{mean-LLN}, we note that from \refeq{IHLLN}, which is now known to be valid for all $n$,
    \eq
    \lbeq{mean-LLNa}
    \Ebold_\beta[\omega_n]=
    -i \nabla \hat{c}_n(0)=-i
    \sum\limits_{j=1}^n \bigl[i\theta_j +\nabla e_j(0)\bigr].
    \en
Since $|e_j(k)|=O(|k|^2)$, we have that $\nabla e_j(0)=0$.
Therefore,
    \eq
    \lbeq{mean-LLNb}
    \Ebold_\beta[\omega_n]=
    \sum\limits_{j=1}^n\theta_j
    =n \theta +\sum\limits_{j=1}^n[\theta_j-\theta].
    \en
By \refeq{theta}, \refeq{thetan} and \refeq{pibdscfirst},
we have that
    \eq
    \lbeq{diffthetas}
    \sum\limits_{j=1}^n[\theta_j-\theta]
    =i\sum_{j=1}^n \sum_{s=j+1}^{\infty} \nabla \hat \pi_s(0)
    =i\sum_{s=2}^\infty (n\wedge (s-1)) \nabla \hat \pi_s(0)
    =O\left(\sum_{s=1}^\infty (n\wedge  s) b_s\right).
    \en

For \refeq{LLN}, let $k \in \R^d$.  Then for $n\ge e^{\delta^{-1}|k|}$ we can apply
\refeq{IHLLN} in the form
    \eq
    \hat c_n\big(kn^{-1}\big)= e^{ikn^{-1}\cdot \theta n} \exp{\Big[\sum\limits_{l=1}^n [ikn^{-1} \cdot
    (\theta_l-\theta)+e_l\big(kn^{-1}\big)]\Big]},
    \en
with
    \eq
    |e_j(kn^{-1})|\leq K\frac{|k|^2}{n^2}\sum_{l=1}^j lb_l.
    \en

By \refeq{diffthetas},
    \eq
    \lbeq{LLNpf1}
    \sum\limits_{l=1}^n ikn^{-1} \cdot
    (\theta_l-\theta)=O\left(\frac{|k|}{n}\sum_{s=1}^\infty (n\wedge s) b_s\right).
    \en
Similarly,
    \eq
    \lbeq{LLNpf2}
    \sum\limits_{j=1}^n|e_j(kn^{-1})|
    \leq K\frac{|k|^2}{n^2} \sum\limits_{j=1}^n \sum_{l=1}^j lb_l
    \leq K\frac{|k|^2}{n^2} \sum_{l=1}^n (n-l+1)lb_l\leq K\frac{|k|^2}{n} \sum_{l=1}^nlb_l.
    \en
Together \refeq{LLNpf1} and \refeq{LLNpf2} prove \refeq{LLN} for $n\ge e^{\delta^{-1}|k|}$.

For $n< e^{\delta^{-1}|k|}$ the result
is trivial by writing
\eq
ikn^{-1}\cdot \omega_n=ik\cdot \theta+O(|k|(L+\theta))=ik\cdot \theta+O(|k|e^{\delta^{-1}|k|}n^{-1}).
\en
\qed

\section{Induction for the central limit theorem}
\label{sec-CLT}
In this section we prove a central limit theorem from the
recurrence relation \refeq{laceexp}, or more precisely its
Fourier transform \refeq{laceexpk}, assuming certain bounds
on the coefficients $\hat{\pi}_m(k)$.  The bounds roughly
correspond to upper bounds on the accuracy of the Taylor
approximation of $\hat{\pi}_m(k)$ up to {\em second} order,
and the argument is an extension of the one in Section \ref{sec-LLN}.
In this section, for a $d\times d$ matrix
$\Sigma$, we define its $L^1$-norm by
    \eq
    \lbeq{matrixnorm}
    |\Sigma|=\sum_{i,j=1}^d|(\Sigma)_{ij}|.
    \en

We start by formulating a general assumption, and prove the main result,
Theorem \ref{thm-1}, under this assumption.

\begin{ASSCLT}
\label{ass:CLT}
There exists a non-increasing sequence $\{b_m\}_{m \ge 1}$
independent of $\beta$ with $b_1\ge 1$, and a constant $\veta$
with $\lim_{\beta\downarrow 0} \veta =0$, such that\\
(i)
    \eq
    \lbeq{pibdscfirstagain}
    \hat{\pi}_m(0)=0, \qquad |\nabla
    \hat{\pi}_m(0)|
    \leq \veta  b_m,
    \qquad
    |\nabla^2\hat{\pi}_m(0)|\leq
    \veta m b_m.
    \en
(ii) for all $k\in [-\pi,\pi]^d$,
    \eq
    \lbeq{pibdscsecond}
    \big| \hat{\pi}_m(k)-k\cdot \nabla
    \hat{\pi}_m(0)\big|\leq
    \veta  |k|^2m b_m,\qquad
    \Big| \hat{\pi}_m(k)-k\cdot \nabla
    \hat{\pi}_m(0)- \frac{1}{2} k^t\nabla^2\hat{\pi}_m(0)k \Big|\leq
    \veta  |k|^3m^2 b_m.
    \en
Moreover, $B^*\equiv\sum_{m=1}^\infty m b_m<\infty$,
and there exists $\gamma\in(0,1/2)$ such that
    \eqalign
    \lbeq{bass}
    &d_n \equiv\sum_{m=2}^{n}m b_m \sum_{l=n+1-m}^{n} b_{l+1}=o(1), \quad\text{as }n \ra \infty, \text{ and}\nn\\
    &a_n\equiv \sum_{m=1}^{n} m^{2+\gamma} b_m\le C_a\sqrt{\frac{n}{\log{(n+1)}}}, \quad \text{for all }n \text{ and some } C_a\ge 1.
    \enalign
\end{ASSCLT}


Similarly to \refeq{Bn-def}, we define
    \eq
    \lbeq{An-def}
    A_n=\sum_{j=1}^n a_j,
    \quad
    D_n=\sum_{j=1}^n d_j,
    \quad
    E_n=\sum_{m=1}^{\infty} (m \wedge n) m b_m.
    \en
We will prove a generalised version of Theorems
\ref{thm-1OERRW} and \ref{thm-1ERW},
which is formulated below:

\begin{theorem}[Central limit theorem]
\label{thm-1}
When Assumption (CLT) holds, there exist $\beta_0 = \beta_0(d)>0$, and
$\theta=\theta(\beta)$, and $\Sigma=\Sigma(\beta)$ such that,
for all $\beta \leq \beta_0$,\\
(a)
\eq
\lbeq{meangen}
    \Ebold_\beta[\omega_n]
    = \theta n \left[1+O\left(\frac{1}{n}\right)\right].
\en
(b)
\eq
\lbeq{Varasygen}
    \Var_\beta(\omega_n)
    = \Sigma n+ O\left(D_n\right)+O\left(E_n\right).
\en
(c) there exists $C>0$ such that for every $k\in \R^d$
    \eq
    \lbeq{CLTgen}
    \log\left(\Ebold_\beta\big[e^{ik\cdot \frac{(\omega_n-\theta n)}{\sqrt{n}}}\big]\right) =
    -\hlf k^t\Sigma k + O(|k|e^{C|k|^2}n^{-1/2})
    +O\big(|k|^3n^{-3/2}A_n\big)
    +O\big(|k|^2 e^{C|k|^2}n^{-1}(D_n+E_n)\big)
    .
    \en
The constants $\theta$ and $\Sigma$ (given by \refeq{theta},
\refeq{sigma}) are model dependent.
\end{theorem}
It is not hard to see that each of the $O$ terms in \refeq{CLTgen}
is indeed an error term when we assume that Assumption (CLT) holds.
However, in the general set-up in Assumption (CLT), it is not
clear to us which term on the right-hand side of \refeq{CLTgen}
is typically the largest.

Note that since $D$ has finite range, there exists a constant $C_2\ge 1$ independent of $\beta$ such that
    \eq
    \lbeq{DTaylor2}
    |\hat{D}(k)-1-ik\cdot \zerospeed-\hlf k^t\zerovar k|\leq C_2|k|^3,
    \en
and let $K_2=2C_2$, which is independent of $\beta$.

Recall \refeq{thetan} and define the following approximation to $\Sigma$:
    \eq
    \Sigma_{n} = \zerovar-\theta_{n}\theta_n^t-\sum_{m=2}^n
    \nabla^2\Big[e^{-}i\theta_n\cdot k(m-1)\hat{\pi}_m(k)\Big]_{k=0}.\lbeq{sigman}
    \en
Let $B\equiv\sum_m b_m$ and $d^*\equiv\sup_n d_n$.

Our induction hypothesis for the central limit theorem
is that the following bound holds for
all $\beta\leq \beta_0$, all $0\leq j\leq n$,
and some $\delta\in(0,1),$ independent of $\beta$:

For $k$ such that $|k|^2\le \delta\log(n\vee 3)/n$,
and some $K$ independent of $\beta$ we can write,
    \eq
    \lbeq{IHCLT}
    \hat c_j(k) = \exp{\Big[\sum\limits_{l=1}^j \big[ik\cdot
    \theta_l-\frac12 k^t\Sigma_l k +r_l(k)\big]\Big]}\quad  \mbox{with}~
    |r_j(k)|\leq K(|k|^2d_{j}+|k|^3a_{j}),
    \en
where again the empty sum appearing when $j=0$ is defined to be zero,
and for $n=0$, \refeq{IHCLT} is assumed to hold for all $k\in [-\pi,\pi]^d$.

The initialisation of the induction ($n=0$ case) holds trivially as $1=e^0$.

{\short{We refer to \cite[Sections 5.1 and 5.2]{HofHol10b} for the advancement of the
induction hypothesis, and the proof of Theorem \protect\ref{thm-1} using the
conclusion of the induction hypothesis.}}

{\extend{\subsection{The CLT induction advanced}
\label{sec-adv}
We follow the same strategy as in Section \ref{sec-LLNadv},
now expanding the Fourier transform one order further.
We fix $n\geq 0$.
The induction step will be achieved as
soon as we are able to write
    \eq
    \frac{\hat{c}_{n+1}(k)}{\hat{c}_n(k)}=
    \exp{\big[ik\cdot \theta_{n+1}-\frac12 k^t\Sigma_{n+1} k+r_{n+1}(k)\big]},
    \en
for $r_{n+1}(k)$ satisfying the required bound. For this,
we write
    \eq \lbeq{eF}
    \frac{\hat{c}_{n+1}(k)}{\hat{c}_{n}(k)}
    =
    1+ik\cdot\theta_{n+1}-\frac12 k^t(\Sigma_{n+1}+\theta_{n+1}\theta_{n+1}^t)k+r'_{n+1}(k)
    \en
and then set
    \eq
    \lbeq{eE}
    r_{n+1}(k) =
    \log{\big[1+ik\cdot \theta_{n+1}-\frac12 k^t(\Sigma_{n+1}+\theta_{n+1}\theta_{n+1}^t)k+r'_{n+1}(k)\big]}
    - ik\cdot \theta_{n+1}+\frac12 k^t\Sigma_{n+1}k.
    \en
The following lemma is an easy consequence of \refeq{zerovardef} and \refeq{sigman}:
\begin{lemma}
\label{lem-sigmabds}
We have $|\zerospeed|\le L$ and $|\zerovar|\le d^2L^2$, and,
when Assumption (CLT) holds, for all $n$, $|\theta_n|\le L+\veta B$, and
    \eqalign
    |\Sigma_n|\le &d^2L^2+(L+\veta B)^2+2d^2(L+\veta B)B^*+\veta B^*, \quad \text{and}\nn\\ |\Sigma_n+\theta_n\theta_n^t|\le &d^2L^2+2d^2(L+\veta B)B^*+\veta B^*.
    \enalign
\end{lemma}

Suppose that the required bound \refeq{IHCLT} holds for $r'_{n+1}(k)$
with constant $K_2$. Then, by the assumption on $a_j$ in \refeq{bass},
we have that, for $k$
satisfying $|k|^2\le \delta \log{(n\vee 3)}/n\leq 2\delta$, and since $\delta<1$,
    \eqn{
    |r'_{n+1}(k)|\leq K_2 \delta (d^*+\sqrt{\delta} C_a)\leq K_2 \delta (d^*+C_a).
    }
Choose $\beta_0$ so that $\veta \le 1$ for all $\beta\le \beta_0$, so that,
by Lemma \ref{lem-sigmabds}, for $k$
satisfying $|k|^2\le \delta \log{(n\vee 3)}/n\leq 2\delta$ in
\refeq{IHCLT}, and using $L,B^*\geq 1$,
    \eqalign
    &|k||\theta_{n+1}|+\frac12 |k|^2|\Sigma_{n+1}+\theta_{n+1}\theta_{n+1}^t|+|r'_{n+1}(k)|\\
    &\qquad\le \sqrt{2\delta}(L+\veta B) +\delta\Big(d^2L^2+(L+\veta B)^2+2d^2(L+\veta B)B^*+\veta B^*\Big)
    +K_2 \delta (d^*+C_a)\nn\\
    &\qquad\leq \sqrt{2\delta}(L+B) +\delta\Big(5d^2(L+\veta B)B^*+K_2(d^*+C_a)\Big)\leq 1/2,\nn
    \enalign
when we take $\delta\leq \delta^*$, which is defined by
    \eq
    \delta^*=\min\Big\{(L+B)^{-2}/32, \left(4\big(5d^2(L+B)B^*+K_2(d^*+C_a)\big)\right)^{-1}\Big\}.
    \en
Therefore we may apply Taylor's Theorem $|\log(1+x)-x+\frac{x^2}{2}|\le 8|x|^3$
for $|x|\le 1/2$ to \refeq{eE}.  This implies
that when the required bound holds for $r'_{n+1}(k)$ with
constant $K_2$, it also holds for $r_{n+1}(k)$ for some $K$
independent of $\beta$, since the terms of order
$k$ and $|k|^2$ in \refeq{eE} cancel.  Specifically, if
$|r'_{n+1}(k)|\le K_2(|k|^2d_{n+1}+|k|^3a_{n+1})$ then using
Taylor's Theorem, followed by the assumed bound on $r'_{n+1}(k)$
, we obtain
\eqalign
|r_{n+1}(k)|\le & |r'_{n+1}(k)|+|k||\theta_{n+1}|\left(\hlf|k|^2|\Sigma_{n+1}+\theta_{n+1}\theta_{n+1}^t|+|r'_{n+1}(k)|\right)\nn\\
& \quad+\hlf\left(\hlf|k|^2|\Sigma_{n+1}+\theta_{n+1}\theta_{n+1}^t|+|r'_{n+1}(k)|\right)^2\nn\\
& \quad+8\left(|k||\theta_{n+1}|+\hlf|k|^2|\Sigma_{n+1}+\theta_{n+1}\theta_{n+1}^t|+|r'_{n+1}(k)|\right)^3\nn\\
\le &CK_2^3(|k|^2d_{n+1}+|k|^3a_{n+1})\le K(|k|^2d_{n+1}+|k|^3a_{n+1}),
\enalign
when $K\ge CK_2^3$.  Here $C\ge 1$ is a constant that depends
on $C_a,B,B^*,d^*,L,d$, but is independent of $\beta$ and $\delta$,
and we have used that $|r'_{n+1}(k)|^2\leq CK_2^2(|k|^2d_{n+1}+|k|^3a_{n+1})$
since $a_{n+1}\ge 1$ and
$|k|^2\le \delta \log{(n\vee 3)}/n$, and similarly for $|r'_{n+1}(k)|^3$.

Most of this section will be devoted to the proof of the following lemma:

\begin{lemma} \label{lem-r'bound} If \refeq{IHCLT} holds for all $j\le n$ and $|k|^{2}\le \delta (n+1)^{-1}\log{((n+1)\vee 3)}$ then for such $k$
   \eqarray
    |r'_{n+1}(k)|  &\leq K_2(|k|^2 d_{n+1}+|k|^3 a_{n+1}).\lbeq{r'bound}
    \enarray
\end{lemma}



\subsubsection{Proof of Lemma \ref{lem-r'bound}}
The proof involves expressing $r'_{n+1}(k)$ as a sum of three terms and showing that each term is bounded in absolute value by the right hand side of \refeq{r'bound}.

Recall \refeq{eF}, then
    \eqarray
    r'_{n+1}(k) & =I+II,\nn
    \lbeq{E'}
    \enarray
where
    \eqarray
    I    & = &[\hat{D}(k)-1-ik\cdot \zerospeed +\frac{1}{2} k^t\zerovar k]
    +\sum_{m=2}^{n+1} \Big[\hat{\pi}_m(k)-
    k\cdot\nabla\hat{\pi}_m(0)-\frac{1}{2}k^t\nabla^2\hat{\pi}_m(0)k\Big],\nn\\
    II   & = &\sum_{m=2}^{n+1}
    \Big[\hat{\pi}_m(k)\big[\frac{\hat{c}_{n+1-m}(k)}{\hat{c}_n(k)}-1\big]
    +k\cdot \nabla\hat{\pi}_m(0)i(m-1)k\cdot \theta_{n+1}\Big].\nn
    \enarray
We will bound $|I|$ and $|II|$, and then choose $\beta_0$ sufficiently small so that $|I|+|II|$ satisfies the bound on the right hand side of \refeq{r'bound}.
By \refeq{DTaylor2} and \refeq{pibdscsecond} in Assumption (CLT), and the fact that $a_{n+1}\ge 1$ we have
    \eqarray |I|    & \leq C_2|k|^3 +\sum_{m=2}^{n+1} \veta  |k|^3m^2 b_m
\le (C_2+\veta )|k|^3a_{n+1}.
    \lbeq{Ifinalbound}
    \enarray

To bound $II$, we first split $II =
II_1+II_2,$ with
    \eqarray
    II_1   & = &\sum_{m=2}^{n+1}
        \left[\hat{\pi}_m(k)
        -k\cdot \nabla\hat{\pi}_m(0)\right]\left[\frac{\hat{c}_{n+1-m}(k)}{\hat{c}_n(k)}-1\right],\nn\\
    II_2   & = &\sum_{m=2}^{n+1}k\cdot\nabla\hat{\pi}_m(0)
    \left[\frac{\hat{c}_{n+1-m}(k)}{\hat{c}_n(k)}-1+i(m-1)k\cdot \theta_{n+1}\right].\nn
    \enarray
For $II_1$, we use the first bound in
\refeq{pibdscsecond} in Assumption (CLT) and
Lemma \ref{lem-Taylor} for $j=0$, to get
    \eqarray \lbeq{III_1est}
    |II_1|   & \leq &\sum_{m=2}^{n+1}
        \left|\hat{\pi}_m(k) - k\cdot \nabla\hat{\pi}_m(0)\right|
        \left|\frac{\hat{c}_{n+1-m}(k)}{\hat{c}_n(k)}-1\right|\nn\\
        & \leq &\sum_{m=2}^{n+1} \veta |k|^2mb_m
        \left|\exp{\left[-\sum_{l=n+2-m}^{n}\Big[ik\cdot \theta_{l}-\frac12 k^t\Sigma_{l}k+r_{l}(k)\Big]\right]}-1
        \right|,\nn\\
    & \leq &\veta |k|^2 \sum_{m=2}^{n+1} mb_m
    e^{\chi_{m,n}(k)}\sum_{l=n+2-m}^{n} \left[|k||\theta_l|+\frac12 |k|^2 |\Sigma_{l}| + |r_l(k)|\right],\nn
    \enarray
with
    \eq
    \lbeq{Chimn}
    \chi_{m,n}(k) =
    \sum_{l=n+2-m}^{n} \Big[\frac{|k|^2}{2}|\Sigma_{l}|+ |r_l(k)|\Big].
    \en
Since $a_n$ is increasing, for $|k|^2\le \delta \log{((n+1)\vee 3)}/(n+1)$,
we have that
    \eq
    \lbeq{Chimn-bd}
    \chi_{m,n}(k)
    \leq m|k|^2(C+Kd^*+a_n|k|)\leq m|k|^2(C+Kd^*+\sqrt{\delta}KC_a),
    \en
where we recall that $d^*=\sup_{n} d_n$. Also,
for $|k|^2\le \delta \log{((n+1)\vee 3)}/(n+1)$,
    \eq
    m|k|^2 \leq \delta \log{(m\vee 3)}
    \frac{\log{((n+1)\vee 3)}}{n+1} \frac{m}{\log{(m\vee 3)}}
    \leq \delta \log{(m\vee 3)},
    \en
since $x\mapsto \frac{\log{(x\vee 3)}}{x}$ is decreasing
for $x\geq 0$. As a result, we obtain that,
with $\nu=\delta(C+Kd^*+\sqrt{\delta}KC_a)$,
    \eq
    \lbeq{Chibd}
    e^{\chi_{m,n}(k)}\leq (m\vee 3)^\nu.
    \en
Note that, by picking $\delta>0$ sufficiently small,
we can make $\nu<\gamma$.

For $|k|^2\le \delta \log{((n+1)\vee 3)}/(n+1)$
it follows from Lemma \ref{lem-sigmabds} and \refeq{IHCLT},
using a similar argument as in \refeq{Chimn-bd}, that
    \eqalign
    |k||\theta_l|+\frac{|k|^2}{2}|\Sigma_{l}| + |r_l(k)|
    \leq \Big((L+\veta B)+(\sqrt{\delta}C+\sqrt{\delta}Kd^*+\delta KC_a)\Big)|k| \equiv C_{\sss K}^{\sss{(1)}}
    |k|\lbeq{expbd}.
    \enalign
Therefore,
    \eqarray
    \lbeq{II1finalbound}
    |II_1|
    & \leq &C_{\sss K}^{\sss{(1)}}\veta |k|^3\sum_{m=2}^{n+1} m^{2+\gamma} b_m
    =C_{\sss K}^{\sss{(1)}}\veta |k|^3 a_{n+1}.
    \enarray

For $II_2$ we use \refeq{IHCLT} and Lemma \ref{lem-Taylor} for $j=1$ to obtain
    \eqarray |II_2|& \leq & |k|\sum_{m=3}^{n+1}\veta  b_m
            \left|\exp{\left[-\sum_{l=n+2-m}^{n}\Big[ik\cdot \theta_{l}-\frac12 k^t\Sigma_{l}k+r_{l}(k)\Big]\right]}
            -1+ik\cdot \theta_{n+1}(m-1)\right|\nn\\
        &\leq &\veta |k|\sum_{m=2}^{n+1}b_m
        \sum_{l=n+2-m}^{n}\Big[|k||\theta_{n+1}-\theta_l|+\frac{|k|^2}{2}|\Sigma_l|+|r_l(k)|\Big]\nn\\
        && \quad +\veta |k|\sum_{m=2}^{n+1} b_m
        \left[\sum_{l=n+2-m}^{n} \Big[ |k||\theta_l|+\frac{|k|^2}{2}|\Sigma_l|+|r_l(k)|\Big]\right]^2e^{\chi_{m,n}(k)}.\nn
        \enarray
For $|k|^2\le \delta \log{((n+1)\vee 3)}/(n+1)$, the second
sum can be bounded, using \refeq{expbd} and \refeq{Chibd}, as
    \eq
    \big(C_{\sss K}^{\sss{(1)}}\big)^23^{\gamma}\veta |k|^3\sum_{m=2}^{n+1} m^{2+\gamma} b_m
    \equiv C_{\sss K}^{\sss{(2)}}\veta |k|^3 a_{n+1}.
    \en
By a similar argument as in \refeq{Chimn-bd},
we have for $|k|^2\le \delta \log{((n+1)\vee 3)}/(n+1),$
    \eq
    \frac{|k|^2}{2}|\Sigma_{l}|+|r_{l}(k)|
    \leq C(1+Kd^*+KC_a \sqrt{\delta})|k|^2
    \equiv C_{\sss K}^{\sss{(3)}}|k|^2.
    \en
Therefore,
    \eq
    \veta |k|\sum_{m=2}^{n+1}b_m
    \sum_{l=n+2-m}^{n}\Big[\frac{|k|^2}{2}|\Sigma_l|+|r_l(k)|\Big]
    \leq C_{\sss K}^{\sss{(3)}}\veta |k|^3 \sum_{m=2}^{n+1}m b_m\le
    C_{\sss K}^{\sss{(3)}}\veta |k|^3 a_{n+1}.
    \en
We continue with the remaining contribution to $II_2$.
Since $\{b_m\}_{m\ge 1}$ is a decreasing sequence,
    \eq
    |\theta_{n+1}-\theta_l|
    \leq \veta \sum_{s=l+1}^{n+1} b_s\leq \veta (n-l+1)b_{l+1}.
    \en
Thus,
    \eq
    \sum_{m=2}^{n+1}b_m\sum_{l=n+2-m}^{n}|\theta_{n+1}-\theta_l|
    \leq \veta\sum_{m=2}^{n+1}m b_m \sum_{l=n+2-m}^{n} b_{l+1}
    =\veta d_{n+1}.
    \en

We conclude that
    \eq
    \lbeq{II2finalbound}
    |II_2|\leq \veta \big(\veta |k|^2 d_{n+1}+(C_{\sss K}^{\sss{(2)}}+C_{\sss K}^{\sss{(3)}})|k|^3 a_{n+1}\big).
    \en
We have shown that
    \eq
    \lbeq{finalI+IIbound}
    |I|+|II|\le
    \veta^2 |k|^2 d_{n+1}+(C_2+\veta+\veta (C_{\sss K}^{\sss{(1)}}+C_{\sss K}^{\sss{(2)}}+C_{\sss K}^{\sss{(3)}}))|k|^3 a_{n+1}.
    \en
Choose $\beta_0$ sufficiently small so that for all $\beta\le \beta_0$,
$\veta+\veta (C_{\sss K}^{\sss{(1)}}+C_{\sss K}^{\sss{(2)}}+C_{\sss K}^{\sss{(3)}})\le \hlf K_2$.
Recall that $K_2=2C_2\geq 1$. Then for $\beta\le \beta_0$,
    \eq
    |I|+|II|\le (C_2+\hlf K_2)(|k|^2 d_{n+1}+|k|^3 a_{n+1})\le K_2(|k|^2 d_{n+1}+|k|^3 a_{n+1}),
     \en
as required.  This completes the proof of Lemma \ref{lem-r'bound}.
\qed

\subsection{Proof of Theorem \protect\ref{thm-1}}
\label{sec-pf1}
We will make use of the following lemma:

\begin{lemma}\label{lem-derE} For all $\beta\leq \beta_0$, and
all $j\in \N$,
\begin{itemize}
\item[(i)] $\nabla r_j(0)=0$, and
\item[(ii)] $|\nabla^2 r_j(0)|  \leq 3Kd^2d_j$.
\end{itemize}
\end{lemma}

\proof The induction hypothesis \refeq{IHCLT}, now verified for
all $j$, states that $|r_j(k)|\leq K(|k|^2d_j+|k|^3a_j)$.
Therefore, letting $[\nabla r_j(0)]_i$ denote the $i^{\rm th}$
coordinate of the vector $\nabla r_j(0)$, we have
    \eq
    \left|[\nabla r_j(0)]_i\right|=\lim_{k_i \rightarrow 0}\frac{|r_j(0,\dots,0,k_i,0,\dots,0)|}{|k_i|}\le \lim_{k_i \rightarrow 0}\frac{K(|k_i|^2d_j+|k_i|^3a_j)}{|k_i|}=0
    \en
Since all partial derivatives of $\hat c_n(k)$ up to second order exist and are continuous, and $\hat c_n(0)=1$, we have from \refeq{eF} and \refeq{eE} that all partial derivatives of $r_j(k)$ up to second order exist in a neighbourhood of $0$ and are continuous.  Let $(\nabla^2 r_j(0))_{lm}$ denote the $(l,m)^{\rm th}$ entry of the matrix $\nabla^2 r_j(0)$ and suppose that $|r_j(k)|\le J_1|k|^2+J_2|k|^3$.  We claim that this implies that $|(\nabla^2 r_j(0))_{lm}|\le 3J_1$ for each $m,l$, from which part {\em (ii)} of the lemma follows immediately.  Without loss of generality we suppose that $l,m \in \{1,2\}$.

Let $h(k_1,k_2)=r_j(k_1,k_2,0,\dots,0)$.  By the second order mean value theorem, $f_{u_1,u_2}(t)\equiv h(tu_1,tu_2)$ satisfies
    \eq
    \lbeq{ft}
    f_{u_1,u_2}(t)=f_{u_1,u_2}(0)+f'_{u_1,u_2}(0)t+f''_{u_1,u_2}(t^*)\frac{t^2}{2}
    \en
for some $t^*\equiv t^*(t,u_1,u_2)\in(0,t)$.

Now $f_{u_1,u_2}(0)=h(0,0)=0$ and
    \eq
    |f'_{u_1,u_2}(0)|=\lim_{t \ra 0}\left|\frac{h(tu_1,tu_2)-h(0,0)}{t}\right|=\lim_{t \ra 0}\left|\frac{h(tu_1,tu_2)}{t}\right|\le \lim_{t \ra 0}\left|\frac{C_{u_1,u_2}(t^2+t^3)}{t}\right|=0
    \en
where we have used the bound on $|r_j(k)|$ in the last inequality.  Thus \refeq{ft} reduces to
    \eq
    f_{u_1,u_2}(t)=f''_{u_1,u_2}(t^*)\frac{t^2}{2},
    \en
and by hypothesis the left hand side is bounded in absolute value by $J_1t^2(u_1^2+u_2^2)+J_2t^3(u_1^2+u_2^2)^{3/2}$.

We now set $t_n=1/n$ and let $t^*_n=t^*(t_n, u_1,u_2)$.  Then for each $n$, $|f''_{u_1,u_2}(t^*_n)|\le 2J_1(u_1^2+u_2^2)+2n^{-1}J_2(u_1^2+u_2^2)^{3/2}$.  By the multivariate chain rule $\frac{d}{dt}h(\vec{g}(t))=\nabla h \cdot \vec{g}'(t)$ we have
    \eq
    f''_{u_1,u_2}(t^*_n)=u_1^2h_{11}(t^*_nu_1,t^*_nu_2)+u_2^2h_{22}(t^*_nu_1,t^*_nu_2)+2u_1u_2h_{12}(t^*_nu_1,t^*_nu_2),
    \en
and thus
    \eq
    \lbeq{partials}
    \big|u_1^2h_{11}(t^*_nu_1,t^*_nu_2)+u_2^2h_{22}(t^*_nu_1,t^*_nu_2)+2u_1u_2h_{12}(t^*_nu_1,t^*_nu_2)\big|\le 2J_1(u_1^2+u_2^2)+\frac{J_2}{n}(u_1^2+u_2^2)^{3/2}.
    \en
Putting $u_1=1,u_2=0$ in \refeq{partials} gives $|h_{11}(t^*_n,0)|\le 2J_1+2n^{-1}J_2$, and similarly $|h_{22}(t^*_n,0)|\le 2J_1+2n^{-1}J_2$.  Letting $n \ra \infty$ and using the fact that $t_n^*\in (0,t_n)$ (so that $t_n=1/n\ra 0$ implies that $t_n^*\ra 0$ as $n\ra \infty$) we have $|h_{11}(0,0)|\le 2J_1$ by continuity of the partial derivatives.  Similarly, by taking $u_1=0,u_2=1$, we obtain $|h_{22}(0,0)|\le 2J_1$.  Next, set $u_1=u_2=1$ in \refeq{partials} and use $|a+b|\le d \Rightarrow |a|\le d+|b|$ to see that
    \eq
    2|h_{12}(t^*_nu_1,t^*_nu_2)|\le |h_{11}(t^*_nu_1,t^*_nu_2)+h_{22}(t^*_nu_1,t^*_nu_2)|
    +2J_1(u_1^2+u_2^2)+\frac{2J_2}{n}(u_1^2+u_2^2)^{3/2}.
    \en
Now use the triangle inequality and let $n\ra \infty$  to get $|h_{12}(0,0)|\le 3J_1$.
\qed

We are now ready to prove the statements in Theorem \ref{thm-1}(a)--(c) one by one.\\
\\
{\bf Proof of Theorem \ref{thm-1}(a):} Using \refeq{IHCLT} and Lemma \ref{lem-derE}(i), we have
    \eqarray
    \sum_{x\in \Z^d} x c_n(x)
    &=&-i \nabla \hat{c}_n(0)=-i
    \sum\limits_{j=1}^n \bigl[i\theta_j +\nabla r_j(0)\bigr]=n\theta+\sum\limits_{j=1}^n [\theta_j-\theta],
    \lbeq{thm-a}
    \enarray
so that it suffices to prove that
    \eq
    \sum\limits_{j=1}^n [\theta_j-\theta]=O(1).
    \en
For this, we use \refeq{theta}, \refeq{thetan} as well as the second
bound in \refeq{pibdscfirstagain} and to note that
    \eq
    \lbeq{thetadiff}
    \sum\limits_{j=1}^n |\theta_j-\theta|\leq
    \sum_{j=1}^n \sum_{m=j+1}^{\infty} |\nabla
    \hat{\pi}_m(0)|
    \leq \veta  \sum_{j=1}^\infty \sum_{m=j+1}^{\infty} b_m
    = \veta  \sum_{m=1}^{\infty} m b_m=O(1),
    \en
by the assumption that $B^*=\sum_{m=1}^{\infty} mb_m<\infty$.
\qed

\medskip \noindent {\bf Proof of Theorem \ref{thm-1}(b):}
Recall that $\mbox{Var}_{\beta}(\omega_n)$ is the covariance
matrix of $\omega_n$.  Then
    \eq
    \lbeq{varlm}
    (\mbox{Var}_{\beta}(\omega_n))_{lm}=\sum_xx_lx_mc_n(x)-\left(\sum_x x_lc_n(x)\right)\left(\sum_x x_m c_n(x)\right).
    \en
By \refeq{IHCLT} and Lemma \ref{lem-derE}(i-ii), and writing $[\theta_p]_l$ for the $l^{\rm th}$ component of $\theta_p$,
    \eqarray
    \sum_{x\in \Z^d} x_lx_m c_n(x)
    &=&-(\nabla^2\hat{c}_n(0))_{lm}\nn\\
    &=&\sum_{p=1}^n\left((\Sigma_p)_{lm}-(\nabla^2r_p(0))_{lm}\right)
    -\sum_{p,q=1}^n\left(i[\theta_p]_l+[\nabla r_p(0)]_l\right)\left(i[\theta_q]_m+[\nabla r_q(0)]_m\right)\nn\\
    &=&\sum_{p=1}^n\left((\Sigma_p)_{lm}+O(d_p)\right)+\sum_{p=1}^n[\theta_p]_l\sum_{q=1}^n[\theta_q]_m.
    \lbeq{secmom}
    \enarray

It follows from \refeq{thm-a} that $\sum_x x_lc_n(x)=\sum_{p=1}^n [\theta_p]_l$ and from \refeq{varlm} and \refeq{secmom} that
    \eqarray
    (\mbox{Var}_{\beta}(\omega_n))_{lm}&=&\sum_{p=1}^n\left((\Sigma_p)_{lm}+O(d_p)\right)\nn\\
&=&n(\Sigma)_{lm}+O\left(\sum_{p=1}^nd_p\right)+\sum_{p=1}^n((\Sigma_p)_{lm}-(\Sigma)_{lm}).\nn
    \enarray
Therefore to complete the proof, it is sufficient to show that for $p \le n$,
    \eq
    |(\Sigma_p)_{lm}-(\Sigma)_{lm}|=O\left(1\vee (p\wedge n)pb_p\right).
    \lbeq{enough}
    \en
By \refeq{sigma} and \refeq{sigman}, the left hand side of \refeq{enough} is bounded by
    \eqarray
    &&|[\theta]_l[\theta]_m-[\theta_p]_l[\theta_p]_m|+\sum_{r=p+1}^{\infty}\left|\left[(\nabla^2e^{-i(r-1)k\cdot \theta} \hat{\pi}_r(k))_{lm}\right]_{k=0}\right|\nn\\
    &&\quad+\sum_{r=2}^p\left|\left[(\nabla^2e^{i(r-1)k\cdot \theta_p} \hat{\pi}_r(k))_{lm}\right]_{k=0}-\left[(\nabla^2e^{i(r-1)k\cdot \theta} \hat{\pi}_r(k))_{lm}\right]_{k=0}\right|\nn\\
    &&\qquad \leq|[\theta]_l||[\theta]_m-[\theta_p]_m|+|[\theta_p]_m||[\theta]_l-[\theta_p]_l|\nn\\
    &&\qquad \quad+\sum_{r=p+1}^{\infty}\left((r-1)\Big(|[\theta]_m||[\nabla \hat{\pi}_r(0)]_l|+|[\theta]_l||[\nabla \hat{\pi}_r(0)]_m|\Big)+|(\nabla^2 \hat{\pi}_r(0))_{lm}|\right)\nn\\
    &&\qquad \quad +\left|\theta_p-\theta\right|\sum_{r=2}^p|[\nabla \hat{\pi}_r(0)]|,
    \enarray
since $\hat{\pi}_r(0)=0$ by \refeq{pibdscfirstagain}.
The first two terms are $O(1)$ using the fact that $|\theta|$ is finite and the $|\theta_p|$ are uniformly bounded together with \refeq{thetadiff}.  By \refeq{pibdscfirstagain}, the third term is bounded by
    \eq
    \veta |\theta|\sum_{p=1}^n\sum_{r=p+1}^{\infty}rb_r\le \veta |\theta|\sum_{r=1}^{\infty}(r\wedge n)rb_r,\nn
    \en
while, again by \refeq{pibdscfirstagain}, the last term is bounded. This completes the proof.\qed

\noindent {\bf Proof of Theorem \protect\ref{thm-1}(c):}
Fix $k \in \R^d$.  Then for $n\ge e^\frac{|k|^2}{\delta}$, we can apply \refeq{IHCLT} in the form
        \eqarray
        \hat{c}_n(n^{-\hlf}k)
        &=&\exp{\Big[\sum_{j=1}^n \big(i n^{-\hlf}k\cdot \theta_j
        -\frac{1}{2}n^{-1}k^t\Sigma_jk + r_j(n^{-\hlf}k)\big)\Big]}\nn\\
        &=&\exp{\big[i
        k\cdot \theta \sqrt{n}
        -\frac{1}{2}k^t\Sigma k\big]}\\
        &&\qquad \qquad \times
        \exp{\big[\sum_{j=1}^n i n^{-\hlf}k \cdot [\theta_j-\theta]
        -\hlf n^{-1}k^t\sum_{j=1}^n [\Sigma_j-\Sigma] k\big]}
        \exp{\sum_{j=1}^n r_j(n^{-\hlf}k)}.\nn
        \lbeq{11cpf}
        \enarray
From \refeq{thetadiff} we have
    \eq
    \Big|\sum_{j=1}^n i n^{-\hlf}k \cdot [\theta_j-\theta]\Big|
    \leq n^{-\hlf}|k| \sum_{j=1}^n |\theta_j-\theta|= O(n^{-\hlf}|k|),
    \en
and using \refeq{enough} we obtain
    \eq
    \Big|\frac{k^t}{2n}\sum_{j=1}^n [\Sigma_j-\Sigma] k\Big|
    \leq \frac{|k|^2}{2n}\sum_{j=1}^n |\Sigma_j-\Sigma|
    =O\Big(\frac{|k|^2}{n} \sum_{m=1}^{\infty} (m \wedge n) m b_m\Big).
    \en
Finally we use \refeq{IHCLT} to get
    \eq
    \sum_{j=1}^n |r_j(n^{-\hlf}k)| \leq O\Big(\frac{|k|^2}{n}\sum_{j=1}^n d_j\Big)
    +O\Big(\frac{|k|^3}{n^{3/2}}\sum_{j=1}^n a_j\Big).
    \en
This proves the bound in Theorem \protect\ref{thm-1}(c) for $n\ge e^\frac{|k|^2}{\delta}$.
The bound holds trivially for $n \le e^\frac{|k|^2}{\delta}$ by writing
\eqarray
ik\cdot \frac{(\omega_n-n\theta)}{\sqrt{n}}+\hlf k^t\Sigma k&=&O(|k|^2+|k|n^{\hlf})=O(|k|^2n^{-1}(D_n+E_n)n+|k|n^{-\hlf}n)\nn\\
&=&O(|k|^2n^{-1}(D_n+E_n)e^{\delta^{-1}|k|^2}+|k|n^{-\hlf}e^{\delta^{-1}|k|^2}).\nn
\enarray

\qed
}}

\section{Bounds on the lace expansion}
\label{sec-pibds}
In this section, we give bounds on the lace expansion coefficients, and
verify that these bounds imply Theorems \ref{thm-1OERRW}, \ref{thm-1ERW}
and \ref{thm-2ERW}. We start in Section \ref{sec-single} by formulating some
general bounds on $\hat \pi_m(0), \nabla \hat \pi_m(0)$ and
$\nabla^2 \hat \pi_m(0)$ that will reduce the bounds on
the derivatives of $\hat \pi_m(k)$ to a single bound{\extend{, which we
will prove separately for each model.}}
{\short{The proofs of these bounds are deferred to \cite[Sections 6.2-6.4]{HofHol10b}.}}
{\extend{In Section \ref{sec-pibdsOERRW}, we prove
the bounds on the lace expansion coefficients for once edge-reinforced
random walk with drift, and complete the proof of Theorem \ref{thm-1OERRW}.
In Section \ref{sec-pibdsERW}, we prove
the bounds on the lace expansion coefficients for excited
random walk, and complete the proof of
Theorems \ref{thm-1ERW}--\ref{thm-2ERW}.  In Section \ref{sec-pibdsRWpRE} we give the corresponding results for the random walk in partially random environment.}}

\subsection{Reduction to a single bound}
\label{sec-single}
Recall \refeq{piNdef} and the definition $\mc{A}_{m,\sN}
=\{(j_1, \dots, j_N)\in \Z_+^N:\sum_{l=1}^Nj_l=m-N-1\}$, and define
   \eqalign
    &\pi_m^{\smallsup{N}}(x,y)
    =&\sum_{\vec{j}\in \mc{A}_{m,N}}\sum_{\walkvec{0}{1}}\sum_{\walkvec{1}{j_{\sss 1}+1}}\dots\sum_{\walkvec{N}{j_{\sss N}+1} }I_{\{\omega^{(N)}_{j_{\sN}}=x, \omega^{(N)}_{j_{\sN}+1}=y\}}D(\walkcoor{0}{1})
    \prod_{n=1}^{N}\Delta^{\sss (n)}_{\sss j_n+1}\prod_{i_{n}=0}^{j_{n}-1}p^{\walkvec{n-1}{j_{n-1}+1}\circ \walkvec{n}{i_{n}}}\left(\walkcoor{n}{i_{n}},\walkcoor{n}{i_{n}+1}\right),  \lbeq{piNxydef}
 \enalign
so that
    \eq
    \pi_m^{\smallsup{N}}(y)=\sum_{x}\pi_m^{\smallsup{N}}(x,y).
    \en
We also let
    \eq
    \pi_m(x,y)=\sum_{N=1}^{\infty} \pi_m^{\smallsup{N}}(x,y).
    \en
The starting point for the bounds on the lace expansion coefficients
for self-interacting random walks is the following proposition:

\begin{PRP}[Reduction of the bounds on the expansion coefficients]
\label{prop-pibdgen}
~\\
For a self-interacting stochastic process with range $L<\infty$,
where $\pi_m^{\smallsup{N}}(y)$ is given by \refeq{piNdef},
the following bounds hold:
    \eqalign
    \hat \pi_m(0)&=0,\lbeq{pibd1}\\
    |\nabla\hat \pi_m(0)|&\leq \sqrt{d}L\sum_{x,y} |\pi_m(x,y)|,
    \lbeq{pibd2}\\
    |\nabla^2\hat \pi_m(0)|&\leq (dL)^2 (2m-1)\sum_{x,y} |\pi_m(x,y)|,
    \lbeq{pibd3}\\
    |\hat{\pi}_m(k)|&\leq |k|L\sum_{x,y}|\pi_m(x,y)|,
    \lbeq{pibd4}\\
    \big|\hat{\pi}_m(k)-k\cdot\nabla\hat{\pi}_m(0)\big|
    &\leq |k|^2 mL^2 \sum_{x,y}|\pi_m(x,y)|,
    \lbeq{pibd5}\\
    \big|\hat{\pi}_m(k)-k\cdot\nabla\hat{\pi}_m(0)-\hlf k\nabla^2\hat \pi_m(0) k^t\big|
    &\leq |k|^3 m^2 L^3 \sum_{x,y}|\pi_m(x,y)|.
    \lbeq{pibd6}
    \enalign
\end{PRP}

\proof
We note that for every $x\in \Z^d$ and $N \ge 1$,
    \eqalign
    \sum_y\Delta^{\smallsup{N}}_{j_{N}+1}I_{\{\omega^{(N)}_{j_N}=x, \omega^{(N)}_{j_N+1}=y\}}&=\sum_y\left(p^{\walkvec{N-1}{j_{N-1}+1}\circ \walkvec{N}{j_N}}(x,y)-p^{\walkvec{N}{j_N}}(x,y)\right)=1-1=0\nn
    \enalign
from which it follows immediately that, for every $x\in \Z^d$,
    \eq
    \lbeq{pixyzero}
    \sum_y\pi_m(x,y)=0.
    \en
Summing \refeq{pixyzero} over $x$ establishes \refeq{pibd1}.
Furthermore, again by \refeq{pixyzero}, we have that
    \eqalign
    \lbeq{1derivpi}
    [\nabla \hat \pi_m(0)]_l&=i\sum_y y_l \pi_m(y)
    =i\sum_{x,y} y_l \pi_m(x,y)=i\sum_{x,y} x_l \pi_m(x,y)
    +i\sum_{x,y} [y_l-x_l] \pi_m(x,y)\nn\\
    &=i\sum_{x,y} [y_l-x_l] \pi_m(x,y).
    \enalign
For walks with range $L$, we have that $|y_j-x_j|\leq L$,
so that
    \eq
    |[\nabla \hat \pi_m(0)]_l|\leq L\sum_{x,y} |\pi_m(x,y)|,
    \en
which establishes \refeq{pibd2} since $\sum_{l=1}^d u_l^2\le d\max_l{|u_l|^2}$.
Similarly,
    \eqalign
    -[\nabla^2 \hat \pi_m(0)]_{st}
    &=\sum_y y_sy_t \pi_m(y)=\sum_{x,y} y_sy_t\pi_m(x,y)\nn\\
    &=\sum_{x,y} x_sx_t \pi_m(x,y)
    +\sum_{x,y} [y_s-x_s]x_t\pi_m(x,y)\nn\\
    &\qquad +\sum_{x,y} [y_t-x_t]x_s\pi_m(x,y)
    +\sum_{x,y} [y_s-x_s][y_t-x_t]\pi_m(x,y)\nn\\
    &=\sum_{x,y} [y_s-x_s]x_t\pi_m(x,y)+\sum_{x,y} [y_t-x_t]x_s\pi_m(x,y)
    +\sum_{x,y} [y_s-x_s][y_t-x_t]\pi_m(x,y).\nn
    \enalign
We use that $|y_j-x_j|\leq L$ and $|x_j|\leq L(m-1)$ to obtain
    \eq
    |[\nabla^2 \hat \pi_m(0)]_{st}|\leq (2m-1)L^2 \sum_{x,y} |\pi_m(x,y)|.\nn
    \en
This establishes \refeq{pibd3} by \refeq{matrixnorm}.


By \refeq{pixyzero},
    \eq
    \lbeq{pik-rewr}
    \hat{\pi}_m(k)=\sum_{x,y} e^{ik\cdot y}\pi_m(x,y)
    =\sum_{x,y} e^{ik\cdot x}[e^{ik\cdot(y-x)}-1]\pi_m(x,y).
    \en
Since $|x-y|\leq L$, this immediately yields \refeq{pibd4}.
Together with \refeq{1derivpi}, \refeq{pik-rewr} gives
    \eqalign
    \hat{\pi}_m(k)-k\cdot\nabla\hat{\pi}_m(0)
    =&\sum_{x,y} \Big[[1+(e^{ik\cdot x}-1)][e^{ik\cdot(y-x)}-1]-ik\cdot(y-x)\Big]\pi_m(x,y),\nn \\
    =&\sum_{x,y}[e^{ik\cdot x}-1][e^{ik\cdot(y-x)}-1]\pi_m(x,y) +\sum_{x,y} [e^{ik\cdot(y-x)}-1-ik\cdot(y-x)]\pi_m(x,y).\nn
    \enalign
By Lemma \ref{lem-Taylor}, $|e^{iu}-1|\le |u|$ and $|e^{iu}-1-iu|\le \hlf u^2$.  Together with the finite range properties of the walk this proves \refeq{pibd5}.
The final claim is proved similarly by first showing that
\eqalign
\hat{\pi}_m(k)-k\cdot\nabla\hat{\pi}_m(0)-\hlf k\nabla^2\hat \pi_m(0)k^t=\sum_{x,y}\Big[& [e^{ik\cdot(y-x)}-1-ik\cdot(y-x)+\hlf k(y-x)^t(y-x)k^t]\nn\\
&+[e^{ik\cdot x}-1-ik\cdot x][e^{ik\cdot(y-x)}-1]\nn\\
&+ik\cdot x[e^{ik\cdot(y-x)}-1-ik\cdot(y-x)]\Big]\pi_m(x,y),\nn
    \enalign
and then using $|e^{iu}-1-iu+\hlf u^2|\le \frac{1}{6}|u|^3$ together with the previous estimates.
\qed
\vskip0.5cm

\noindent
We conclude that the bounds in \refeq{pibdscfirst},
\refeq{piLLNbound}, \refeq{pibdscfirstagain} and \refeq{pibdscsecond}
follow if we can show that
    \eq
    \lbeq{aimpibds}
    \sum_{x,y} |\pi_m(x,y)| \leq \veta  b_m,
    \en
for some sequence $\{b_m\}_{m \ge 1}$ satisfying
the appropriate conditions formulated in Assumptions (LLN) and (CLT).
{\ch{In the following proposition, we state the precise form of our
bounds on the lace expansion coefficients, followed by the proof
of our main results subject to these bounds.}
\ch{
\begin{PRP}[Bounds on the expansion coefficients for each of our models]
\label{prop-pibds-general}
~\\
(a) For OERRWD, there exist $\beta_0>0$ and $\mc{J}>0$ such that for all $|\beta|\le \beta_0$,
    \eq
    \lbeq{piNbound1}
    \sum_{x,y} |\pi_m(x,y)|\le C\beta e^{-\mc{J}m},
    \en
where $\mc{J}$ depends on $\beta_0, w_0, d$ but is independent of $\beta$.\\
(b) For ERW with {\bf $d-1>4$}, there exist $\beta_0>0$ and $C>0$ such that for all $\beta\le \beta_0$
    \eq
    \sum_{x,y}|\pi_m(x,y)|\le \frac{C\beta}{(m+1)^{\frac{d-3}{2}}}.
    \en
(c) For RWpRE with $d_1>4$, there exist $\beta_0>0$ and $C>0$ such that for all $\beta\le \beta_0$
    \eq
    \sum_{x,y}|\pi_m(x,y)|\le \frac{C\beta}{(m+1)^{\frac{d_1-2}{2}}}.
    \en
\end{PRP}
We now complete the proofs of our main results subject to Proposition \ref{prop-pibds-general}:}
\medskip

\noindent
\ch{{\it Proof of Theorem \ref{thm-1OERRW} subject to
Proposition \ref{prop-pibds-general}(a).} We use Proposition \ref{prop-pibdgen}
and \ref{prop-pibds-general} as well as Theorem \ref{thm-1} to complete the
proof of Theorem \ref{thm-1OERRW}. When $b_m=e^{-\mc{J} m}$, Assumption (CLT) is
satisfied.  Also, \refeq{mean} is directly implied by
\refeq{meangen}. Furthermore, the error terms in \refeq{Varasygen} can all be seen
to be $O(1)$, which proves \refeq{var}. Finally, for each
$k$, as $n\ra \infty$  \refeq{CLTgen} implies that
    \eq
    \Ebold_\beta[e^{ik\cdot (\omega_n-\theta n)/\sqrt{n}}] \ra e^{-\hlf k^t\Sigma k}.
    \en
Clearly, this implies \refeq{CLT}.
\qed}
\medskip

\noindent
\ch{{\em Proof of Theorems \ref{thm-1ERW} and \ref{thm-2ERW} subject to
Proposition  \ref{prop-pibds-general}(b)}.
By Propositions \ref{prop-pibdgen} and \ref{prop-pibds-general}(b),
Assumption (LLN) holds with $\veta =C\beta$ and
$b_m=(m+1)^{-(d-3)/2}$ when $d>5$ (i.e.~$d-1>4$).  Thus, Theorem
\ref{thm-LLNgen} applies, and it is an easy exercise to
see that when $b_m=(m+1)^{-(d-3)/2}$ and $d>5$,
the error terms given in \refeq{mean-LLN} and
\refeq{LLN} are sufficient to prove Theorem \ref{thm-2ERW}.

Similarly, Propositions \ref{prop-pibdgen} and \ref{prop-pibds-general}(b)
show that Assumption (CLT) holds with $\veta =C\beta$ and
$b_m=(m+1)^{-(d-3)/2}$ when $d>7$ (i.e.~$d-1>6$).  Thus, Theorem \ref{thm-1}
applies and we now show that when
$b_m=(m+1)^{-(d-3)/2}$ and $d>8$, the error terms given in
\refeq{Varasygen} and \refeq{CLTgen} are sufficient to prove
Theorem \ref{thm-1ERW}.

Indeed, note that
$E_n=\sum_{m=1}^\infty (m\wedge n)m b_m=O(n^{-(d-9)/2}\log{n})=o(n)$,
when $d>7$. Furthermore, by \cite[Lemma 3.2]{HofSla02} and the fact that
$(d-3)/2>1$ when $d>5$, we obtain that $d_n=O(n^{-(d-5)/2})$,
so that $D_n=\sum_{m=1}^n d_m=O(1)$ for $d>7$. Finally,
$a_n=\sum_{m=1}^n m^{2+\gamma} b_m=O(n^{\gamma-(d-9)/2}\vee 1)$,
which is $O(n^{c})$ for some $c<1/2$ when $d>8$ (i.e.~$d-1>7$) and $\gamma$ is sufficiently small. In this case,
also $A_n=\sum_{m=1}^n a_m =o(n^{3/2})$. This identifies all error
terms in \refeq{Varasygen} and \refeq{CLTgen}.
\qed}
\medskip

\noindent
\ch{{\it Proof of Theorems \ref{thm-1RWpRE} and \ref{thm-2RWpRE} subject to
Proposition \ref{prop-pibds-general}(c).} Theorems \ref{thm-1RWpRE} and
\ref{thm-2RWpRE} follow exactly as in the proofs of Theorems \ref{thm-1ERW}
and \ref{thm-2ERW}, when $d_1>7$ and $d_1>4$, respectively.\\
\qed}

{\short{We refer to \cite[Sections 6.2, 6.3 and 6.4]{HofHol10b}
for the proofs of Proposition \ref{prop-pibds-general}(a), (b) and (c), respectively.
In \cite[Section 6.5]{HofHol10b}, these bounds are discussed in more detail.}}

{\extend{We will prove \refeq{aimpibds} for once-edge-reinforced random walk
in Section \ref{sec-pibdsOERRW} and for excited random walk in
Section \ref{sec-pibdsERW} below.}}

{\extend{
\subsection{Bounds for once-edge-reinforced random walk}
\label{sec-pibdsOERRW}
In this section we prove Proposition \ref{prop-pibds-general}(a).
The bounds in this section are based on the following large deviations estimates.
\begin{lemma}[Large deviations]
\label{lem:ldev}
Whenever $\zerospeed\ne 0$, there exist $\beta_0=\beta_0(D(\cdot),w_0(\cdot))>0$ and $\mc{I}=\mc{I}(D(\cdot),w_0(\cdot))>0$ such that for all $|\beta|\le \beta_0$,
\eqalign
\lbeq{ld1}&\sup_{\vec{\eta}}\Qbold^{\vec{\eta}}_{\beta}(\omega_n=\omega_0)\le e^{-\mc{I}n}, \quad \text{and}\\
\lbeq{ld2}&\sup_{z,\walkvec{i-2}{j_{i-2}+1}}
\sum_{\walkvec{i-1}{j_{i-1}+1}}\mQ_{\beta}^{\walkvec{i-2}{j_{i-2}+1}}(\vec{\omega}_{j_{i-1}}=\walkvec{i-1}{j_{i-1}})I_{\{\walkcoor{i-1}{j_{i-1}}=z\}}
\mQ_{\beta}^{\walkvec{i-1}{j_{i-1}+1}}(\omega_{j_i}=\walkcoor{i-1}{l})\le Ke^{-\mc{I}(j_{i-1}-l+j_i)},
\enalign
where the supremum is over all $(j_{i-2}+1)$-step random walk paths
$\walkvec{i-2}{j_{i-2}+1}$, and $K$ is a constant that depends only on $L,d$.  The law of the $i^{\rm th}$ walk $\walkvec{i}{j_{i}+1}$ depends on
the $(i-1)^{\rm st}$ walk $\walkvec{i-1}{j_{i-1}+1}$ but not on the $(i-2)^{\rm nd}$ walk $\walkvec{i-2}{j_{i-2}+1}$.
\end{lemma}
\proof Under $\mathbb{Q}_0$, $\vec{\omega}$ is a simple random walk with bounded increments and non-zero drift (without loss of generality assume the drift is in the positive coordinate direction(s)).  It follows that for $z$ sufficiently small (and negative)
    \eqn{
    \mathbb{E}_{\mathbb{Q}_0}[\exp\{z\cdot(\omega_1-\omega_0)\}]
    =1+z\cdot\mathbb{E}_{\mathbb{Q}_0}[\omega_1-\omega_0]+O(z^2L^2)<1.
    }
Thus, by Cram\'er's Theorem (e.g.\ see  \cite[Theorem 2.2.30]{DemZei98}) there exists $J=J(D(\cdot),w_0(\cdot))>0$ such that $\mathbb{Q}_{0}^{}(\omega_n=\omega_0)\le e^{-Jn}$ for all $n$.
Let $\Omega$ denote the support of $D(x)$.  It is easy to show that for every $\beta\in [0, \beta_0]$ and $\vec{\nu}$,
    \eq
    \lbeq{transbound}
    p^{\vec{\nu}}(x,y)\le \left(1+ C\beta_0\right)D(y-x)
    \en
when $C \ge 1/w_0(x,y)$ (similarly for $\beta\in [-\beta_0,0]$ when $C\ge (|\Omega|-1)/(\sum_{u \sim x}w_0(x,u) -\beta_0(|\Omega|-1))$.

By translation invariance, $w_0(\cdot)\ge W$ is uniformly bounded
from below as a function on $\Omega$.  We fix
    \eqalign
    C\ge & \max\left\{\frac{|\Omega|-1}{\hlf\sum_{u \sim 0}w_0(0,u)}, \sup_{y\sim 0}\frac{1}{w_0(0,y)}\right\}, \quad  \text{and}\nn\\
    \beta_0\le &\min\left\{\frac{\sum_{u\sim 0}w_0(0,u)}{2(|\Omega|-1)},J/(2C)\right\},
    \enalign
where the constant $C>0$ shall be determined in the course of the proof, and recall that
    \eq
    \mathbb{Q}_{\beta}^{\vec{\eta}}(
    \omega_n=x)=\sum_{\walk_{n}:\omega_n=x}\prod_{i=0}^{n-1}p^{\walk_i\circ \vec{\eta}}(\walk_{i+1}-\walk_{i}).
    \en
The bound \refeq{ld1} with $\mc{I}=J/2$ follows immediately from this by \refeq{transbound} by choosing $\beta_0$ sufficiently small so that $\log(1+C\beta_0)\le J/2$.

The second bound is obtained similarly, using
\refeq{transbound} after the $l^{\rm th}$ step of $\walkvec{i-1}{}$, with the constant arising from the ``missing" transition probability corresponding to the sum over $\walkcoor{i-1}{j_{i-1}+1}$.
This proves \refeq{ld2} with $\mc{I}=J/2$.
\qed


{\noindent \em Proof of Proposition \ref{prop-pibds-general}(a).} We bound $\sum_{x,y} |\pi^{\sss (N)}_m(x,y)|$ and sum the resulting bound over $N$.  For $N=1$, $m\ge 2$, \refeq{piNxydef} and \refeq{DeltaERRW} give
   \eqalign
    \sum_{x,y} |\pi^{\sss (1)}_m(x,y)|&\le    \sum_{x,y}\sum_{\walkcoor{0}{1}}D(\walkcoor{0}{1})\sum_{\walkvec{1}{m-1}}\Qbold_{\beta}^{\walkvec{0}{1}}(\vec{\omega}_{m-2}=\walkvec{1}{m-2})|\Delta^{\smallsup{1}}_{m-1}|I_{\{\omega^{(1)}_{m-2}=x\}}I_{\{\omega^{(1)}_{m-1}=y\}}\nn \\
&\le C\beta \mathbb{Q}_{\beta}(\omega_{m-2}=\omega_0)\le C\beta e^{-\mc{I}(m-2)}\le C\beta e^{-\mc{I}m},
    \enalign
where we have applied the first bound of Lemma \ref{lem:ldev} in the last line, and the value of $C$ changes from place to place.

For general $N$, we have that
\eqalign
    \lbeq{sumpiN}
    \sum_{x,y} |\pi^{\sss (N)}_m(x,y)|&\le
     &\sum_{\vec{j}\in \mc{A}_{m,N}}\sum_{\walkvec{0}{1}}\sum_{\walkvec{1}{j_{\sss 1}+1}}\dots\sum_{\walkvec{N}{j_{\sss N}+1} }D(\walkcoor{0}{1})
    \prod_{n=1}^{N}|\Delta^{\sss (n)}_{\sss j_n+1}|\mQ_{\beta}^{\walkvec{n-1}{j_{n-1}+1}}(\vec{\omega}_{j_n}=\walkvec{n}{j_{n}}),
    \enalign
where, by \refeq{DeltaERRW},
    \eq
    \lbeq{Deltabound}
    \sum_{\walkcoor{i}{j_i+1}}|\Delta^{\smallsup{i}}_{j_{i}+1}|\le C\beta \sum_{l_{i-1}=0}^{j_{i-1}}I_{\{\walkcoor{i}{j_i}=\walkcoor{i-1}{l_{i-1}}\}}.
    \en

Let $N\ge 2$, and for $q\in \{0,1\}$ let $A_q=\{i\le N:(N-i)\mod 2=q\}$ and $B_q$
be the set of $\vec{j}\in \mc{A}_{m,\sN}$ such that $\sum_{i \in A_q}(j_i+1)\ge m/2$.
For $r=0, \dots, N-1$, denote by $l_r\le j_r$ the number of
steps in the $r^{\rm th}$ walk $\walkvec{r}{j_{r}+1}$
up to the intersection point as in \refeq{Deltabound} (in particular, $l_0=0$).
Then, combining \refeq{sumpiN} and \refeq{Deltabound},
\eqalign
\lbeq{sumpiN2}
\sum_{x,y} |\pi^{\sss (N)}_m(x,y)|\le &(C\beta)^N\sum_{\vec{j}\in \mc{A}_{m,N}}\sum_{\vec{l}}
\sum_{\walkvec{0}{1}}\sum_{\walkvec{1}{j_{\sss 1}+1}}\dots\sum_{\walkvec{N}{j_{\sss N}+1} }D(\walkcoor{0}{1})
    \prod_{n=1}^{N}I_{\{\walkcoor{n}{j_n}=\walkcoor{n-1}{l_{n-1}}\}}\mQ_{\beta}^{\walkvec{n-1}{j_{n-1}+1}}(\vec{\omega}_{j_n}=\walkvec{n}{j_{n}}).
\enalign
The bound is now split into four cases, depending on whether $N$ is even or odd,
and on whether $\vec{j}\in B_0$ or $\vec{j}\in B_1\setminus B_0$.  See Figure \ref{fig-4cases}.

\begin{figure}
\includegraphics[scale=.75]{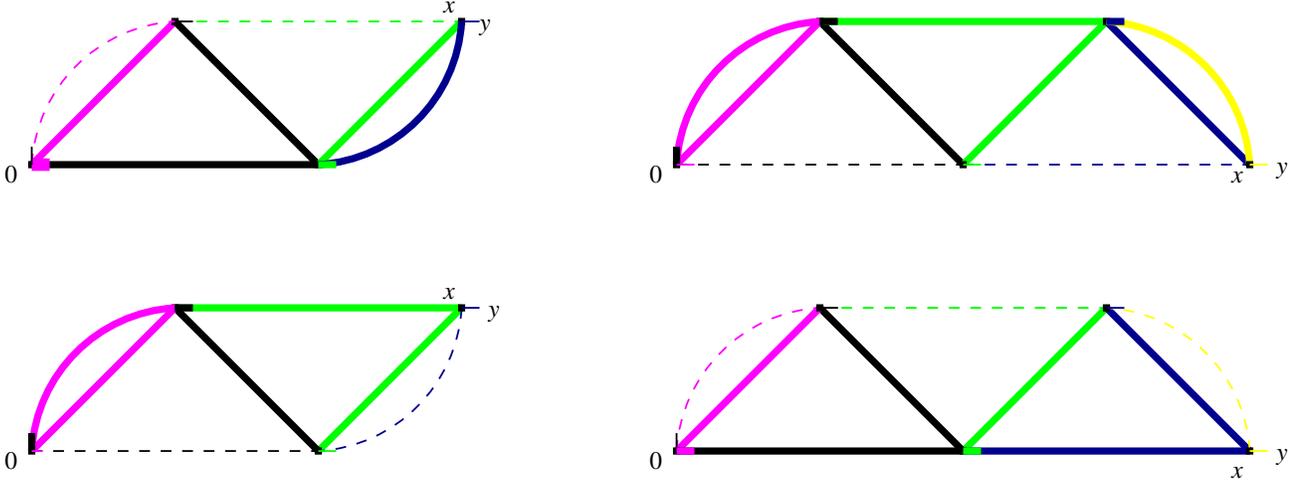}
\caption{Illustration of the four cases of the diagrammatic bounds for OERRWD.  On the left $N$ is even ($N=4$) and on the right $N$ is odd ($N=5$).  In the first row $\vec{j}\in B_0$, while in the second, $\vec{j}\in B_1\setminus B_0$. In each case the thick lines indicate the loops (whose total length is of order half of the total length $m$ of the diagram) that give an exponentially small bound.}
\label{fig-4cases}
\end{figure}

{\bf (a) The bound for $N$ even and $\vec{j}\in B_0$.}
When $N$ is even we bound the contribution to \refeq{sumpiN} from $\vec{j}\in B_0$
by using the following two bounds, the first of which follows immediately from the second bound of Lemma \ref{lem:ldev}, while the last holds (with equality) trivially.

The first fact is that for each even $i \in [2, N]$, uniformly in $\walkvec{i-2}{j_{i-2}+1}$,
    \eqalign
    \lbeq{NevenB0i}
    \sum_{\walkvec{i-1}{j_{i-1}+1} }I_{\{\walkcoor{i-1}{j_{i-1}}=\walkcoor{i-2}{l_{i-2}}\}}\mQ_{\beta}^{\walkvec{i-2}{j_{i-2}+1}}(\vec{\omega}_{j_{i-1}}=\walkvec{i-1}{j_{i-1}})
   \sum_{\walkvec{i}{j_{i}+1} }I_{\{\walkcoor{i}{j_i}=\walkcoor{i-1}{l_{i-1}}\}}\mQ_{\beta}^{\walkvec{i-1}{j_{i-1}+1}}(\vec{\omega}_{j_i}=\walkvec{i}{j_{i}})
    \le Ce^{-\mc{I}(j_i+(j_{i-1}-l_{i-1}))}.
    \enalign
The second fact is that,
    \eqalign
    \lbeq{NevenB01}
    \sum_{\walkvec{0}{1}}D(\walkcoor{0}{1})\le 1.
    \enalign

By successive applications of \refeq{NevenB0i} and
lastly \refeq{NevenB01}, when $N$ is even we obtain a bound
on the contribution to \refeq{sumpiN} from $\vec{j}\in B_0$ (whence $\sum_{i\le N, \rm{even}}j_i\ge (m-N)/2$), of
    \eqalign
    \lbeq{NevenB0aa}
    (C\beta)^N\sum_{\vec{j}\in B_0}\sum_{\vec{l}}\prod_{2\le i\le N,~ \rm{even}}e^{-\mc{I}(j_i+(j_{i-1}-l_{i-1}))}&\le (C\beta)^Ne^{-\mc{I}m/2}\sum_{\vec{j}\in B_0}\sum_{\vec{l}}\prod_{2\le i\le N,~ \rm{even}}e^{-\mc{I}(j_{i-1}-l_{i-1})},
    \nn\\
    \enalign
where the constant has changed (to accommodate a factor $e^{\mc{I}N/2}$).  Using the fact that there are at most $j_i+1$ possible values $\{0,1,\dots,j_{i}\}$ for $l_i$, this is bounded above by
    \eqalign
    \lbeq{NevenB0ab}
    &(C\beta)^Ne^{-\mc{I}m/2}\sum_{\vec{j}}\prod_{i=1}^{N}(j_i+1),
    \enalign
which in turn can be bounded by the integral
    \eqalign
    \lbeq{NevenB0ac}
    &e^{-\mc{I}m/2}(C\beta)^N\int_{0}^{m+3}x_1\int_{0}^{m+3-x_1}x_2\cdots
    \int_{0}^{m+3-(x_1+\dots+x_{N-1})}x_{\sN}dx_{\sN}\cdots dx_1.
    \enalign
It is an easy exercise in integration by parts that
    \eq
    \lbeq{intbyparts}
    \int_0^{a-\sum_{i=1}^{j-1}x_i}\frac{x_{j}}{(2(N-j))!}\left(a-\sum_{i=1}^{j}x_i\right)^{2(N-j)}dx_{j}
    =\frac{1}{(2(N-(j-1)))!}\left(a-\sum_{i=1}^{j-1}x_i\right)^{2(N-(j-1))}.
    \en
Applying \refeq{intbyparts} $N$ times, we bound \refeq{NevenB0ac} by
    \eqalign
    \lbeq{evenandoddbound}
    &  e^{-\mc{I}m/2}(C\beta)^N\frac{(m+3)^{2N}}{(2N)!}\le e^{-\mc{I}m/2}C\beta(C\beta)^{N/2}\frac{(m+3)^{2N}}{(2N)!}.
    \enalign
{\bf (b) The bound for $N$ even and $\vec{j}\in B_1\setminus B_0$.}
When $N$ is even we bound the contribution to \refeq{sumpiN} from $\vec{j}\in B_1\setminus B_0$ by using the following three facts, the first of which is obtained by simply evaluating the sum, while the second and third follow immediately from Lemma \ref{lem:ldev}.

The first fact is that uniformly in $\walkvec{N-1}{j_{N-1}+1}$,
\eqalign
\lbeq{NevenB1N}
&\sum_{\walkvec{N}{j_{\sss N}+1} }I_{\{\walkcoor{N}{j_N}=\walkcoor{N-1}{l_{N-1}}\}}\mQ_{\beta}^{\walkvec{N-1}{j_{N-1}+1}}(\vec{\omega}_{j_N}=\walkvec{N}{j_{N}})\le C,
   \enalign
where the constant (which depends only on $L,d$) is a result of summing over $\sum_{\walkcoor{N}{j_{\sss N}+1} }$.
The second fact is that for each odd $i \in [3, N-1]$, uniformly in $\walkvec{i-2}{j_{i-2}+1}$, \refeq{NevenB0i} holds. The third fact is that
    \eqalign
    \lbeq{NevenB11}
    \sum_{\walkvec{0}{1}}D(\walkcoor{0}{1})
    \Qbold^{\walkvec{0}{1}}(\walkcoor{1}{j_1}=\walkcoor{0}{l_{0}})
    \le e^{-\mc{I}j_1}
    \enalign

By first applying \refeq{NevenB1N}, followed by successive applications of \refeq{NevenB0i} and lastly \refeq{NevenB11}, when $N$ is even we obtain a bound on the contribution to \refeq{sumpiN} from $\vec{j}\in B_1$, of
    \eqalign
    \lbeq{NevenB1}
    (C\beta)^N\prod_{1\le i\le N-1,~ \rm{odd}}\sum_{\vec{l}}e^{-\mc{I}(j_i+(j_{i-1}-l_{i-1}))}&\le
    (C\beta)^Ne^{-\mc{I}m/2}\sum_{\vec{j}\in B_1}\sum_{\vec{l}}
    \prod_{1\le i\le N-1,~ odd}e^{-\mc{I}(j_{i-1}-l_{i-1})}\nn\\
    &\le (C\beta)^Ne^{-\mc{I}m/2}\sum_{\vec{j}}\prod_{i=1}^N(j_i+1),
    \enalign
which is bounded by \refeq{evenandoddbound} just as in the previous case.
\\
{\bf (c),(d) The bounds for $N$ odd.} The bounds for $N\geq 3$ odd are similar to the
bounds described above, and we will omit the details.
When $N$ is odd, we bound the contribution from $\vec{j}\in B_0$
by using the bound \refeq{NevenB0i} successively for each $i \in [3,N]$,
and finally \refeq{NevenB11}.  For $ \vec{j} \in B_1\setminus B_0$,
we use \refeq{NevenB1N}, \refeq{NevenB0i} and finally \refeq{NevenB01}.
In both cases we obtain the same bound \refeq{evenandoddbound}.
\\

To complete the proof of Proposition \ref{prop-pibds-general}(a),
we sum \refeq{evenandoddbound} over $N\geq 2$, giving at most
    \eq
    \lbeq{2ndlast}
    C\beta e^{-(\mc{I}/2-(C\beta)^{1/4})m}.\\
    \en
Choosing $\beta_0$ sufficiently small so that $(C\beta)^{1/4}\le \mc{I}/4$ for all $|\beta|\le \beta_0$, we have that \refeq{2ndlast} is bounded by $C\beta e^{-\mc{J}m}$, where $\mc{J}=\mc{I}/4$ is independent of $\beta$.
\qed
}}

{\extend{
\subsection{Bounds for excited random walk}
\label{sec-pibdsERW}
In this section we prove Proposition \ref{prop-pibds-general}(b).

In bounding the diagrams arising from the expansion applied to excited random walk, we will make use of the following lemma, in which $\mathbb{Q}^{\vec{\eta}}$ denotes the law of an excited random walk with history $\vec{\eta}$, where $\vec{\eta}$ is a finite path:
\begin{lemma}
\label{lem:d-1}
For excited random walk in $d>2$ dimensions,
    \eq
    \sup_{x,\vec{\eta}}\mathbb{Q}^{\vec{\eta}}(\omega_m=x)\le \frac{C}{(m+1)^{\frac{d-1}{2}}}.
    \en
\end{lemma}
\proof
Let $Y_m=\#\{k\le m:\omega_k\notin \{\omega_{k-1}\pm e_1\}\}$ denote the number of steps taken in the dimensions $2,\dots, d$ by the excited random walk up to time $m$.  Note that for excited random walk and simple random walk, $Y_n$ has the same distribution.  Then $Y_m \sim Bin(m,q)$ where $q=(d-1)/d>\hlf$ for $d>2$, and standard large deviations estimates give $\mathbb{P}(Y_m<m/2)\le e^{-mI}$ for some $I>0$.

Now for each $\vec{\eta}$, with endpoint $u$,
    \eqalign
    \lbeq{d-1m}
    \mathbb{Q}^{\vec{\eta}}(\omega_m=x)\le & \mathbb{Q}^{\vec{\eta}}(\omega_m^{[2,\dots,d]}=x^{[2,\dots,d]})=\mathbb{P}_u(\omega_m^{[2,\dots,d]}=x^{[2,\dots,d]}),
    \enalign
where $\mathbb{P}_u$ denotes the law of a simple random walk starting at $u$.  For $m$ even, this is bounded by
    \eqalign
    \lbeq{d-1meven}
    \mathbb{P}_0(\omega_m^{[2,\dots,d]}=0^{[2,\dots,d]})\le  & \sum_{k=m/2}^m\mathbb{P}_0(\omega_m^{[2,\dots,d]}=0^{[2,\dots,d]}|Y_m=k)\mathbb{P}(Y_m=k)+\mathbb{P}(Y_m<m/2)\nn\\
    \le &\sum_{k=m/2}^m \frac{C}{(k+1)^{\frac{d-1}{2}}}\mathbb{P}(Y_m=k)+e^{-Im}\nn \\
    \le &\frac{C}{(m+1)^{\frac{d-1}{2}}}\sum_{k=m/2}^m \mathbb{P}(Y_m=k)+e^{-Im}\le \frac{C}{(m+1)^{\frac{d-1}{2}}}.
    \enalign
For $m$ odd, \refeq{d-1m} is bounded by $2d\mathbb{P}_0(\omega_{m+1}^{[2,\dots,d]}=0^{[2,\dots,d]})$ and we proceed as in \refeq{d-1meven}.\qed

\vspace{.5cm}
Recall that $\mc{A}_{m,\sN}\equiv \{\vec{j}\in \Z_+^{N}:\sum j_i =m-N-1\}$, and that for $N\ge 1$,
    \eqalign
    \lbeq{erwpia}
\sum_{x,y} |\pi^{\sss (N)}_m(x,y)|\le &(C\beta)^N\sup_{\vec{\eta}}\sum_{\vec{j}\in \mc{A}_{m,N}}\sum_{\vec{l}}
\sum_{\walkvec{0}{1}}\sum_{\walkvec{1}{j_{\sss 1}+1}}\dots\sum_{\walkvec{N}{j_{\sss N}+1} }p^{\sss \vec{\eta}}(u,\walkcoor{0}{1})
    \prod_{n=1}^{N}I_{\{\walkcoor{n}{j_n}=\walkcoor{n-1}{l_{n-1}}\}}
    \mQ_{\beta}^{\walkvec{n-1}{j_{n-1}+1}}(\vec{\omega}_{j_n}=\walkvec{n}{j_{n}})\nn\\
 =&(C\beta)^N\sup_{\vec{\eta}} \Pi_m^{\smallsup{N},\vec{\eta}},
\enalign
where $u$ is the endpoint of the finite path $\vec{\eta}$, and we take this expression as the definition of $\Pi_m^{\smallsup{N},\vec{\eta}}$.  See the top diagram in Figure \ref{fig-bigPi}.
\begin{PRP}[Bounds on the expansion coefficients for ERW]
\label{prp:erwpi}
For excited random walk with $d>5$, the following bound holds:
    \eq
    \sum_{x,y}|\pi_m^{\smallsup{N}}(x,y)|\le \frac{(C\beta)^N}{(m+1)^{\frac{d-3}{2}}}.
    \en
\end{PRP}

\begin{figure}
\includegraphics[scale=.6]{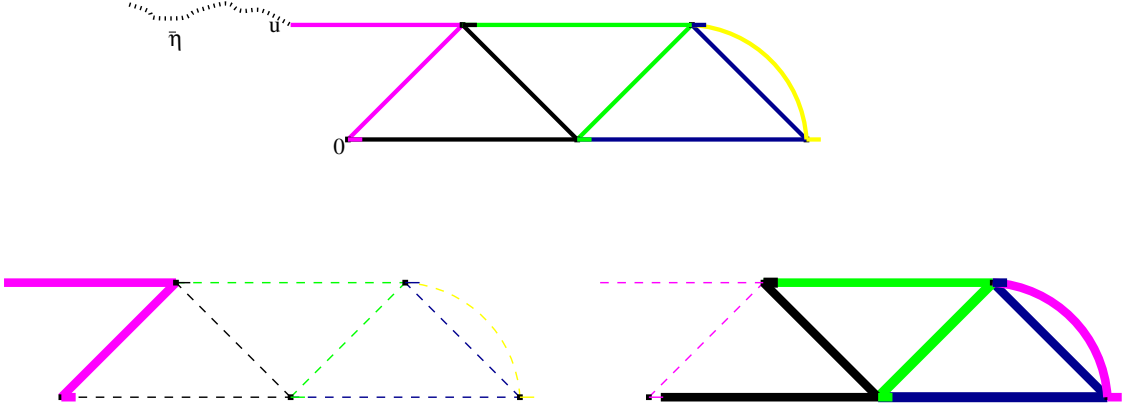}
\caption{A diagrammatic representation of $\Pi_m^{\smallsup{N},\vec{\eta}}(u)$ for $N=5$, followed by the decomposition of the diagram when $j_1>m/2$ and when $j_1\le m/2$ respectively.  In each case, the induction hypothesis is applied to the subdiagram of length $m-(j_1+1)$ that excludes the first walk, and the required decay comes from the part of the diagram with thick lines.}
\label{fig-bigPi}
\end{figure}

In view of \refeq{erwpia}, the conclusion of Proposition \ref{prp:erwpi} follows immediately from the following lemma:
\begin{lemma}
\label{lem:bigPI}
For $d>5$, there exists $C$ independent of $\beta$ such that
    \eq
  \lbeq{lembigpi}  \sup_{u,\vec{\eta}}\Pi_m^{\smallsup{N},\vec{\eta}}\le \frac{C^N}{(m+1)^{\frac{d-3}{2}}}.
    \en
\end{lemma}
\proof We first prove by induction on $N\geq 1$ that
\eq
\lbeq{pieces1}
 \sum_{\vec{j}\in \mc{A}_{m,N}}\sum_{l_1,\dots,l_{N-1}}
\sum_{\walkvec{1}{j_{\sss 1}+1}}\dots\sum_{\walkvec{N}{j_{\sss N}+1} }
    \prod_{n=1}^{N}I_{\{\walkcoor{n}{j_n}=\walkcoor{n-1}{l_{n-1}}\}}\mQ_{\beta}^{\walkvec{n-1}{j_{n-1}+1}}(\vec{\omega}_{j_n}=\walkvec{n}{j_{n}}) \le \frac{C^N}{(m+1)^{\frac{d-3}{2}}}.
    \en
For $N=1$, $j_1=m-2$ and \refeq{pieces1} is less than or equal to
\eq
\lbeq{pieces2}
\sup_{v,\vec{\eta}}\sum_{\walkvec{1}{m-1}}I_{\{\walkcoor{1}{m-2}=v\}}\mQ_{\beta}^{\vec{\eta}}(\vec{\omega}_{m-2}=\walkvec{1}{m-2})=C\sup_{v,\vec{\eta}}\mQ_{\beta}^{\vec{\eta}}(\omega_{m-2}=v)
 \le \frac{C}{(m+1)^{\frac{d-1}{2}}},
\en
where the first constant arises from the sum over $\walkcoor{1}{m-1}$.

For $N\ge 2$, \refeq{pieces1} is bounded by
\eqalign
&\sup_{v,\vec{\eta}}\sum_{j_1\le m-2}\sum_{l_1\le j_1}\sum_{\walkvec{1}{j_{\sss 1}+1}}I_{\{\walkcoor{1}{j_1}=v\}}\mQ_{\beta}^{\vec{\eta}}(\vec{\omega}_{j_1}=\walkvec{1}{j_{1}})\nn\\
&\times \sum_{\vec{j}'\in \mc{A}_{m-j_1-1,N-1}}\sum_{l_2,\dots,l_{N-1}}
\sum_{\walkvec{2}{j_{\sss 2}+1}}\dots\sum_{\walkvec{N}{j_{\sss N}+1} }
    \prod_{n=2}^{N}I_{\{\walkcoor{n}{j_n}=\walkcoor{n-1}{l_{n-1}}\}}\mQ_{\beta}^{\walkvec{n-1}{j_{n-1}+1}}(\vec{\omega}_{j_n}=\walkvec{n}{j_{n}})\nn\\
    \le &\sup_{v,\vec{\eta}}\sum_{j_1\le m-2}\frac{C^{N-1}}{(m-j_1-1)^{\frac{d-3}{2}}}\sum_{l_1\le j_1}\sum_{\walkvec{1}{j_{\sss 1}+1}}I_{\{\walkcoor{1}{j_1}=v\}}\mQ_{\beta}^{\vec{\eta}}(\vec{\omega}_{j_1}=\walkvec{1}{j_{1}})\lbeq{pieces3}\\
    \le &\sum_{j_1\le m-2}\frac{C^{N-1}}{(m-j_1)^{\frac{d-3}{2}}}j_1\frac{C'}{(j_1+1)^{\frac{d-1}{2}}}\le C^{N-1}\sum_{j_1\le m-2}\frac{1}{(m-j_1)^{\frac{d-3}{2}}}\frac{C'}{(j_1+1)^{\frac{d-3}{2}}},\lbeq{pieces4}
    \enalign
using the induction hypothesis to get \refeq{pieces3} and \refeq{pieces2} to get \refeq{pieces4}.  The result \refeq{pieces1} now follows by splitting the sum over $j_1$ into the cases $j_1\le m/2+1$ and $j_1>m/2+1$, taking the term of order $m^{-(d-3)/2}$ outside the sum and performing the remaining finite (since $d>5$) sum.

To prove \refeq{lembigpi}, for $N=1$, we have from \refeq{erwpia} and Lemma \ref{lem:d-1} that
    \eqalign
    \Pi_m^{\smallsup{1},\vec{\eta}}=&
\sum_{\walkvec{0}{1}}\sum_{\walkvec{1}{m-2}}p^{\sss \vec{\eta}}(u,\walkvec{0}{1})
    I_{\{\walkcoor{1}{m-3}=\walkcoor{0}{0}\}}\mQ_{\beta}^{\walkvec{0}{1}}(\vec{\omega}_{m-3}=\walkvec{1}{m-3})\nn\\
\le & \sup_{v,\vec{\eta}'}\mQ_{\beta}^{\vec{\eta}}(\omega_{m-3}=v)\sum_{\walkcoor{1}{m-2}}1\le \frac{C}{(m+1)^\frac{d-1}{2}},\lbeq{pieces5}
    \enalign
uniformly in $\vec{\eta}$, which initialises the induction hypothesis.

For $N\ge 2$, and for any $\vec{\eta}$, proceeding as in the proof of \refeq{pieces1},
    \eqalign
    \lbeq{bigPiinduct}
    \Pi_m^{\smallsup{N},\vec{\eta}}\le &
    \sum_{j_1\le m-2}\sum_{l_1\le j_1}\sum_{\walkvec{0}{1}}\sum_{\walkvec{1}{j_{\sss 1}+1}}p^{\sss \vec{\eta}}(u,\walkvec{0}{1})
    I_{\{\walkcoor{1}{j_1}=\walkcoor{0}{0}\}}\mQ_{\beta}^{\walkvec{0}{1}}(\vec{\omega}_{j_1}=\walkvec{1}{j_{1}})\nn\\
 &\times   \sum_{\vec{j}'\in \mc{A}_{m-(j_1+1),N-1}}\sum_{\vec{l}'}
\sum_{\walkvec{2}{j_{\sss 2}+1}}\dots\sum_{\walkvec{N}{j_{\sss N}+1} }
    \prod_{n=2}^{N}I_{\{\walkcoor{n}{j_n}=\walkcoor{n-1}{l_{n-1}}\}}\mQ_{\beta}^{\walkvec{n-1}{j_{n-1}+1}}(\vec{\omega}_{j_n}=\walkvec{n}{j_{n}})\\
    \le &C^{N-1}\sum_{j_1\le m-N-1}\frac{1}{(m-j_1)^{\frac{d-3}{2}}}j_1\frac{C'}{(j_1+1)^\frac{d-1}{2}},
    \enalign
    using \refeq{pieces5} and \refeq{pieces1}.  The result follows as for \refeq{pieces4}.
This completes the proof of Lemma \ref{lem:bigPI}, and hence also Proposition \ref{prp:erwpi}.
\qed
}}

{\extend{
\subsection{Bounds for random walk in a partially random environment}
\label{sec-pibdsRWpRE}
In this section we prove Proposition \ref{prop-pibds-general}(c), proceeding similarly to the excited random walk case.
The main ingredient needed is the analogue of Lemma \ref{lem:d-1} for RWpRE, which is the following Lemma.
\begin{lemma}
\label{lem:d_1}
For RWpRE with $d_1\ge 1$ dimensions,
    \eq
    \sup_{x,\vec{\eta}}\mathbb{Q}^{\vec{\eta}}(\omega_m=x)\le \frac{C}{(m+1)^{\frac{d_1}{2}}}.
    \en
\end{lemma}
\proof
Let $Y_m=\#\{k\le m:\omega_k\notin \{\omega_{k-1}\pm e_i, i=1,\dots, d_0\}\}$ denote the number of steps taken in the dimensions $d_0+1,\dots d$ by the RWpRE up to time $m$.  Then there exists a sequence of random variables $Y_m'\sim Bin(m,\delta)$ such that $Y_m'\le Y_m$ for all $m$, and $\mP(Y_m<m\frac{\delta}{2})\le \mP(Y_m'<m\frac{\delta}{2})\le e^{-Im}$ for some $I>0$, by standard large deviations estimates.  Now proceed as in the proof of Lemma \ref{lem:d-1} to get the result.\qed

Since for RWpRE as defined in Section \ref{sec-RWRE}, the $\Delta$ factors satisfy the same bounds  (\ref{e:DeltaRWRE}) as excited random walk (\ref{e:DeltaExRW}), the analysis continues exactly as in Section \ref{sec-pibdsERW} except that the exponents have changed in Proposition \ref{prp:erwpi} and Lemma \ref{lem:bigPI} from $\frac{d-3}{2}=\frac{d-1}{2}-1$ to $\frac{d_1}{2} -1$.  In the inductive analysis we then use the fact that when $\frac{d_1}{2}-1>1$  (i.e.~$d_1>4$),
\[\sum_{j\le m-2}\frac{1}{(m-j)^{\frac{d_1}{2} -1}} j \frac{1}{(j+1)^{\frac{d_1}{2}}}\le \ch{\frac{C}{(m+1)^{\frac{d_1}{2}-1}}}.\]

}}

{\extend{
\subsection{Discussion of the bounds}
\label{sec-disc-bds}
In the examples given in this paper, an estimate of the form
    \eq
    \lbeq{missingbound}
    \sup_{\vec{\eta},x}\mathbb{Q}^{\vec{\eta}}(\omega_m=x)\le A(m)
    \en
is crucially used in bounding the diagrams, where $A(m)$ is
decreasing sufficiently rapidly in $m$.  In the case of the
reinforced random walk with drift, Cram\'er's Theorem enabled
such a result with $A(m)$ exponentially small in $m$.  For
excited random walk, the simple random walk behaviour in
all but the first dimension gave such a result with
$A(m)=(m+1)^{-(d-1)/2}$.  Similarly for random walk in partially random environment with $A(m)=(m+1)^{-d_1/2}$.  In these examples, we
ignore considerable information contained in the expansion
in order to bound certain quantities arising from the expansion
in terms of diagrams.  In the case of excited random walk, we
bounded these diagrams using very simple, but non-optimal
estimates.  The diagrammatic estimates are used to verify
a set of non-optimal assumptions under which the central
limit theorem holds.  Improvements in any of these areas
could lead to a reduction in the dimension above which
our methods imply a central limit theorem for excited
random walk. We note that different bounds,
valid for {\it all} $\beta\in [0,1]$, are proved in \cite{HHmono08}
for ERW in order to prove monotonicity of $\beta\mapsto \theta(\beta,d)$
when $d\geq 9$.


The approach taken above works more generally.  We can obtain a LLN and CLT for any translation invariant self-interacting random walk model that has the properties that
\begin{itemize}
\item[(1)] $p^{\vec{\eta}_m\circ\vec{x}_n}(x_n,x_{n+1})-p^{\vec{x}_n}(x_n,x_{n+1})\ne 0 \Rightarrow x_n\in \vec{\eta}_m$,
\item[(2)] this difference in transition probabilities is small (uniformly) for all possible histories, and
\item[(2)] the walker is ``sufficiently transient" (uniformly) for all possible histories,
\end{itemize}
can be handled in the same way as we have handled the models above.  For an explicit example, one can take an (annealed) multi-cookie random walk in an i.i.d.~random cookie environment with multi-dimensional excitement, provided that there are $d_1>4$ (sufficiently transient) coordinates where the walker is behaving as a simple random walk.

It would require a great advance in our understanding and analysis of the recursion equation, in order for us to apply this methodology to a ``non-repulsive" model such as the once reinforced random walk.
Inductive arguments as in
\cite{HofHolSla98, HofSla02} have been used
rather successfully for oriented percolation \cite{HofSla03b},
the contact process \cite{HofSak04}, and various related problems.  However all of these made crucial use of the self-repellent nature of the problems involved.

}}

{\extend{\section{Proof of the variance formula in Theorem \ref{thm-var-form}}
\label{sec-thm-var-form}
Multiplying both sides of \refeq{laceexp} by $x^{[i]}x^{[j]}=(x^{[i]}-y^{[i]}+y^{[i]})(x^{[j]}-y^{[j]}+y^{[j]})$ and summing over $x$ we obtain,

\eqalign
\lbeq{e1}
\mE[\omega^{[i]}_{n+1}\omega^{[j]}_{n+1}]=&\sum_{y}y^{[i]}y^{[j]}D(y)\sum_xc_n(x-y)+\sum_{y}y^{[i]}D(y)\sum_x(x^{[j]}-y^{[j]})c_n(x-y)\nn\\
&+\sum_{y}y^{[j]}D(y)\sum_x(x^{[i]}-y^{[i]})c_n(x-y)+\sum_{y}D(y)\sum_x(x^{[i]}-y^{[i]})(x^{[j]}-y^{[j]})c_n(x-y)\nn\\
&+\sum_{m=2}^{n+1}\sum_{y}y^{[i]}y^{[j]}\pi_m(y)\sum_xc_{n+1-m}(x-y)+\sum_{m=2}^{n+1}\sum_{y}y^{[i]}\pi_m(y)\sum_x(x^{[j]}-y^{[j]})c_{n+1-m}(x-y)\nn\\
&+\sum_{m=2}^{n+1}\sum_{y}y^{[j]}\pi_m(y)\sum_x(x^{[i]}-y^{[i]})c_{n+1-m}(x-y)\nn\\
&+\sum_{m=2}^{n+1}\sum_{y}\pi_m(y)\sum_x(x^{[i]}-y^{[i]})(x^{[j]}-y^{[j]})c_{n+1-m}(x-y).
\enalign
Since $\sum_x c_n(x)=1=\sum_yD(y)$ and $\sum_y\pi_m(y)=0$, many terms simplify, so that \refeq{e1} becomes
\eqalign
\lbeq{e2}
\mE[\omega^{[i]}_{n+1}\omega^{[j]}_{n+1}]=&\sum_{y}y^{[i]}y^{[j]}D(y)+\sum_{y}y^{[i]}D(y)\sum_x x^{[j]}c_n(x)+\sum_{y}y^{[j]}D(y)\sum_x x^{[i]}c_n(x)\nn\\
&+\sum_x x^{[i]}x^{[j]}c_n(x)+\sum_{m=2}^{n+1}\sum_{y}y^{[i]}y^{[j]}\pi_m(y)+\sum_{m=2}^{n+1}\sum_{y}y^{[i]}\pi_m(y)\sum_x x^{[j]}c_{n+1-m}(x)\nn\\
&+\sum_{m=2}^{n+1}\sum_{y}y^{[j]}\pi_m(y)\sum_x x^{[i]}c_{n+1-m}(x)\nn\\
=&\mE[\omega^{[i]}_{1}\omega^{[j]}_{1}]+\mE[\omega^{[i]}_{1}]\mE[\omega^{[j]}_{n}]+\mE[\omega^{[j]}_{1}]\mE[\omega^{[i]}_{n}]+\mE[\omega^{[i]}_{n}\omega^{[j]}_{n}]\nn\\
&+\sum_{m=2}^{n+1}\sum_{y}y^{[i]}y^{[j]}\pi_m(y)+\sum_{m=2}^{n+1}a_m^{[i]}\mE[\omega^{[j]}_{n+1-m}]+\sum_{m=2}^{n+1}a_m^{[j]}\mE[\omega^{[i]}_{n+1-m}].
\enalign

Turning this into a statement about covariances we have, with
$C(X,Y)=\mE[XY]-\mE[X]\mE[Y]$ denoting the covariance between the random variables $X$ and $Y$,
\eqalign
\lbeq{e3}
C(\omega^{[i]}_{n+1},\omega^{[j]}_{n+1})-C(\omega^{[i]}_{n},\omega^{[j]}_{n})\nn
=&\mE[\omega^{[i]}_{1}\omega^{[j]}_{1}]+\sum_{m=2}^{n+1}\sum_{y}y^{[i]}y^{[j]}\pi_m(y)\\
&+\mE[\omega^{[i]}_{1}]\mE[\omega^{[j]}_{n}]+\mE[\omega^{[j]}_{1}]\mE[\omega^{[i]}_{n}]-\mE[\omega^{[i]}_{n+1}]\mE[\omega^{[j]}_{n+1}]+\mE[\omega^{[i]}_{n}]\mE[\omega^{[j]}_{n}]\nn\\
&+\sum_{m=2}^{n+1}a_m^{[i]}\mE[\omega^{[j]}_{n+1-m}]+\sum_{m=2}^{n+1}a_m^{[j]}\mE[\omega^{[i]}_{n+1-m}].
\enalign
If the right hand side converges then by the first condition of \refeq{lim_exists} it must converge to $\Sigma_{ij}$, since the left hand side summed from $n=0$ to $k-1$ is $C(\omega^{[i]}_{k},\omega^{[j]}_{k})$.  Note that if $\mE[\omega_{n}]=0$ for each $n$ then the last two lines of \refeq{e3} are zero and the claimed result then follows immediately (with $\theta^{[i]}=0$ for all $i$).  Otherwise we need to show that under the conditions of \refeq{lim_exists} and (ii) the right hand side of \refeq{e3} converges to that of \refeq{lim_covariance}.

Use the relationship
\[\mE[\omega_{n+1}^{[i]}]=\mE[\omega_{1}^{[i]}]+\sum_{m=2}^{n+1}a_m^{[i]}+\mE[\omega_{n}^{[i]}]\equiv \theta_{n+1}^{[i]}+\mE[\omega_{n}^{[i]}],\]
to see that
\eqalign
\mE[\omega^{[i]}_{n+1}]\mE[\omega^{[j]}_{n+1}]=&\left(\theta_{n+1}^{[i]}+\mE[\omega_{n}^{[i]}]\right)\left(\theta_{n+1}^{[j]}+\mE[\omega_{n}^{[j]}]\right)\nn\\
=&\theta_{n+1}^{[i]}\theta_{n+1}^{[j]}+\theta_{n+1}^{[i]}\mE[\omega_{n}^{[j]}]+\theta_{n+1}^{[j]}\mE[\omega_{n}^{[i]}]+\mE[\omega_{n}^{[i]}]\mE[\omega_{n}^{[j]}].
\enalign
Thus the right hand side of \refeq{e3} is
\eqalign
\lbeq{e4}
&\mE[\omega^{[i]}_{1}\omega^{[j]}_{1}]+\sum_{m=2}^{n+1}\sum_{y}y^{[i]}y^{[j]}\pi_m(y)+\sum_{m=2}^{n+1}a_m^{[i]}\mE[\omega^{[j]}_{n+1-m}]+\sum_{m=2}^{n+1}a_m^{[j]}\mE[\omega^{[i]}_{n+1-m}]\nn\\
&\quad+\mE[\omega^{[i]}_{1}]\mE[\omega^{[j]}_{n}]+\mE[\omega^{[j]}_{1}]\mE[\omega^{[i]}_{n}]-\theta_{n+1}^{[i]}\theta_{n+1}^{[j]}-\theta_{n+1}^{[i]}\mE[\omega_{n}^{[j]}]-\theta_{n+1}^{[j]}\mE[\omega_{n}^{[i]}]\nn\\
&=\mE[\omega^{[i]}_{1}\omega^{[j]}_{1}]+\sum_{m=2}^{n+1}\sum_{y}y^{[i]}y^{[j]}\pi_m(y)+\sum_{m=2}^{n+1}a_m^{[i]}\mE[\omega^{[j]}_{n+1-m}]+\sum_{m=2}^{n+1}a_m^{[j]}\mE[\omega^{[i]}_{n+1-m}]\nn\\
&\quad-\mE[\omega_{n}^{[j]}]\sum_{m=2}^{n+1}a_m^{[i]}-\mE[\omega_{n}^{[i]}]\sum_{m=2}^{n+1}a_m^{[j]}-\theta_{n+1}^{[i]}\theta_{n+1}^{[j]}.
\enalign
Collecting terms, we can rewrite \refeq{e3} as
\eqalign
C(\omega^{[i]}_{n+1},\omega^{[j]}_{n+1})-C(\omega^{[i]}_{n},\omega^{[j]}_{n})=&\mE[\omega^{[i]}_{1}\omega^{[j]}_{1}]-\theta_{n+1}^{[i]}\theta_{n+1}^{[j]}+\sum_{m=2}^{n+1}a_m^{[i]}\left(\mE[\omega^{[j]}_{n+1-m}]-\mE[\omega_{n}^{[j]}]\right)\nn\\
&+\sum_{m=2}^{n+1}\sum_{y}y^{[i]}y^{[j]}\pi_m(y)+\sum_{m=2}^{n+1}a_m^{[j]}\left(\mE[\omega^{[i]}_{n+1-m}]-\mE[\omega_{n}^{[i]}]\right)\nn\\
=&\mE[\omega^{[i]}_{1}\omega^{[j]}_{1}]-\theta_{n+1}^{[i]}\theta_{n+1}^{[j]}-\sum_{m=2}^{n+1}a_m^{[i]}\sum_{r=n+2-m}^n\theta_r^{[j]}\nn\\
&+\sum_{m=2}^{n+1}\sum_{y}y^{[i]}y^{[j]}\pi_m(y)-\sum_{m=2}^{n+1}a_m^{[j]}\sum_{r=n+2-m}^n\theta_r^{[i]}.\nn
\enalign
The right hand side is equal to
\eqalign
\lbeq{e6}
&\mE[\omega^{[i]}_{1}\omega^{[j]}_{1}]-\theta_{n+1}^{[i]}\theta_{n+1}^{[j]}-\sum_{m=2}^{n+1}a_m^{[i]}\sum_{r=n+2-m}^n\left(\theta^{[j]}-\sum_{k=r+1}^{\infty}a_m^{[j]}\right)\nn\\
&\quad+\sum_{m=2}^{n+1}\sum_{y}y^{[i]}y^{[j]}\pi_m(y)-\sum_{m=2}^{n+1}a_m^{[j]}\sum_{r=n+2-m}^n\left(\theta^{[i]}-\sum_{k=r+1}^{\infty}a_k^{[i]}\right)\nn\\
&=\mE[\omega^{[i]}_{1}\omega^{[j]}_{1}]-\theta_{n+1}^{[i]}\theta_{n+1}^{[j]}-\theta^{[j]}\sum_{m=2}^{n+1}a_m^{[i]}(m-1)+\sum_{m=2}^{n+1}\sum_{y}y^{[i]}y^{[j]}\pi_m(y)-\theta^{[i]}\sum_{m=2}^{n+1}a_m^{[j]}(m-1)\\
&\quad+\sum_{m=2}^{n+1}a_m^{[i]}\sum_{r=n+2-m}^n\sum_{k=r+1}^{\infty}a_k^{[j]}+\sum_{m=2}^{n+1}a_m^{[j]}\sum_{r=n+2-m}^n\sum_{k=r+1}^{\infty}a_k^{[j]}.\nn
\enalign

The first line of the last equality of \refeq{e6} converges to \refeq{lim_covariance}.
It therefore remains to show that the terms on the second line of the last equality of \refeq{e6} converge to zero.  Since $i$ and $j$ are arbitrary, it suffices to verify the result for the first term on the second line of the last equality of \refeq{e6}.  For $n\ge 4$ this term is equal to
\eqalign
\lbeq{zeroterms1}
\sum_{m=2}^{n+1}a_m^{[i]}\sum_{r=n+2-m}^n\sum_{k=r+1}^{\infty}a_k^{[j]}=&\sum_{k=2}^{\infty}a_k^{[j]}\sum_{m=(n+3-k)\vee 2}^{n+1}a_m^{[i]}((k-1)+(m-n-1)).
\enalign
This is bounded in absolute value by
\eqalign
\lbeq{zeroterms2}
&\sum_{k=2}^{\infty}|a_k^{[j]}|(k-1)\sum_{m=(n+3-k)\vee 2}^{n+1}|a_m^{[i]}| +\sum_{k=2}^{\infty}|a_k^{[j]}|\sum_{m=(n+3-k)\vee 2}^{n+1}|a_m^{[i]}|(m-1)
\enalign
The first term of \refeq{zeroterms2} is
\eqalign
&\sum_{k=2}^{\floor{n/2}}|a_k^{[j]}|(k-1)\sum_{m=(n+3-k)\vee 2}^{n+1}|a_m^{[i]}| +\sum_{k=\floor{n/2}+1}^{\infty}|a_k^{[j]}|(k-1)\sum_{m=(n+3-k)\vee 2}^{n+1}|a_m^{[i]}|\nn\\
&\quad \le \sum_{k=2}^{\infty}|a_k^{[j]}|(k-1)\sum_{m=n/2}^{n+1}|a_m^{[i]}| +\sum_{k=\floor{n/2}+1}^{\infty}|a_k^{[j]}|(k-1)\sum_{m=2}^{\infty}|a_m^{[i]}|,
\enalign
which converges to $0$ as $n\ra \infty$, since each of these is the tail of a convergent series multiplied by a convergent series.  Similarly the second term of \refeq{zeroterms1} converges to 0.
\qed}}

\paragraph{Acknowledgements.}
The work of RvdH and MH was supported in part by Netherlands Organisation for
Scientific Research (NWO).  The work of MH was performed in
part at Eindhoven University of Technology.
RvdH thanks Vlada Limic for various discussions
and encouragements at the start of this project.

\def\cprime{$'$}


\begin{thebibliography}{10}

\bibitem{AntRed05}
T.~Antal and S.~Redner.
\newblock The excited random walk in one dimension.
\newblock {\em J. Phys. A}, {\bf 38}(12):2555--2577, (2005).

\bibitem{BS08}
A.-L.~Basdevant and A.~Singh.
\newblock On the speed of a cookie random walk.
\newblock {\em Probab. Theory Related Fields}, {\bf 141}:62--645, (2008).

\bibitem{BenWil03}
I.~Benjamini and D.B. Wilson.
\newblock Excited random walk.
\newblock {\em Electron. Comm. Probab.}, {\bf 8}:86--92 , (2003).


\bibitem{BR07}
J.~B\'{e}rard and A.F.~Ram\'{i}rez.
\newblock Central limit theorem for the excited random walk in dimensions $d\ge2$.
\newblock {\em Electron. Comm. Probab.}, {\bf 12}:303--314 (2007).

\bibitem{BSZ03}
E.~Bolthausen and A.-S.~Sznitman and O.~Zeitouni,
\newblock Cut points and diffusive random walks in random environment.
\newblock {\em Ann. Inst. H. Poincar\'e Probab. Statist.}, {\bf 39}(3):527--555 (2003).

{\extend{\bibitem{BrySpe85}
D.C. Brydges and T.~Spencer.
\newblock Self-avoiding walk in 5 or more dimensions.
\newblock {\em Commun. Math. Phys.}, {\bf 97}:125--148, (1985).
}}

\bibitem{Davi99}
B.~Davis.
\newblock Brownian motion and random walk perturbed at extrema.
\newblock {\em Probab. Theory Related Fields}, {\bf 113}(4):501--518, (1999).

\bibitem{DemZei98}
A.~Dembo and O.~Zeitouni.
\newblock {\em Large deviations techniques and applications}, volume~{\bf 38}
  of {\em Applications of Mathematics (New York)}.
\newblock Springer-Verlag, New York, second edition, (1998).

{\extend{\bibitem{DerSla97}
E.~Derbez and G.~Slade.
\newblock Lattice trees and super-{Brownian} motion.
\newblock {\em Canad.\ Math.\ Bull.}, {\bf 40}:19--38, (1997).
}}

{\extend{\bibitem{DerSla98}
E.~Derbez and G.~Slade.
\newblock The scaling limit of lattice trees in high dimensions.
\newblock {\em Commun.\ Math.\ Phys.}, {\bf 193}:69--104, (1998).
}}

\bibitem{DurKesLim02}
R.~Durrett, H.~Kesten, and V.~Limic.
\newblock Once edge-reinforced random walk on a tree.
\newblock {\em Probab. Theory Related Fields}, {\bf 122}(4):567--592, (2002).

{\extend{\bibitem{HarSla90a}
T.~Hara and G.~Slade.
\newblock Mean-field critical behaviour for percolation in high dimensions.
\newblock {\em Commun. Math. Phys.}, {\bf 128}:333--391, (1990).
}}

{\extend{\bibitem{HarSla90b}
T.~Hara and G.~Slade.
\newblock On the upper critical dimension of lattice trees and lattice animals.
\newblock {\em J. Stat. Phys.}, {\bf 59}:1469--1510, (1990).
}}

\bibitem{HarSla92b}
T.~Hara and G.~Slade.
\newblock The lace expansion for self-avoiding walk in five or more dimensions.
\newblock {\em Reviews in Math.\ Phys.}, {\bf 4}:235--327, (1992).

{\extend{\bibitem{HarSla00a}
T.~Hara and G.~Slade.
\newblock The scaling limit of the incipient infinite cluster in
  high-dimensional percolation. {I}. {C}ritical exponents.
\newblock {\em J. Statist. Phys.}, {\bf 99}(5-6):1075--1168, (2000).
}}

{\extend{\bibitem{HarSla00b}
T.~Hara and G.~Slade.
\newblock The scaling limit of the incipient infinite cluster in
  high-dimensional percolation. {II}. {I}ntegrated super-{B}rownian excursion.
\newblock {\em J. Math. Phys.}, {\bf 41}(3):1244--1293, (2000).
}}

\bibitem{Hofs01}
R.~van~der Hofstad.
\newblock The lace expansion approach to ballistic behaviour for
  one-dimensional weakly self-avoiding walks.
\newblock {\em Probab. Theory Related Fields}, {\bf 119}(3):311--349, (2001).

\bibitem{HofHolSla98}
R.~van~der Hofstad, F.~den Hollander, and G.~Slade.
\newblock A new inductive approach to the lace expansion for self-avoiding
  walks.
\newblock {\em Probab. Theory Related Fields}, {\bf 111}(2):253--286, (1998).

\bibitem{HHmono08}
R.~van~der Hofstad and M.~Holmes.
\newblock Monotonicity for excited random walk in high dimensions.
\newblock Preprint (2008). \ch{To appear in {\em Probab. Theory Related Fields}.}

{\extend{\bibitem{HofHol10a}
R.~van~der Hofstad and M.~Holmes.
\newblock An expansion for self-interacting random walks.
\newblock Submitted. Preprint (2010).}}

{\short{\bibitem{HofHol10b}
R.~van~der Hofstad and M.~Holmes.
\newblock An expansion for self-interacting random walks: extended version.
\newblock (2010).}}

{\extend{\bibitem{HofSak04}
R.~van~der Hofstad and A.~Sakai.
\newblock Gaussian scaling for the critical spread-out contact process above
  the upper critical dimension.
\newblock {\em Electron. J. Probab.}, {\bf 9}:710--769 (electronic), (2004).
}}

{\extend{\bibitem{HofSak08}
R.~van~der Hofstad and A.~Sakai.
\newblock   Convergence of the critical finite-range contact process to super-Brownian motion above the upper critical dimension: I. The higher-point functions.
\newblock Preprint (2009).
}}

\bibitem{HofSla02}
R.~van~der Hofstad and G.~Slade.
\newblock A generalised inductive approach to the lace expansion.
\newblock {\em Probab. Theory Related Fields}, {\bf 122}(3):389--430, (2002).

{\extend{\bibitem{HofSla03b}
R.~van~der Hofstad and G.~Slade.
\newblock Convergence of critical oriented percolation to super-{B}rownian
  motion above {$4+1$} dimensions.
\newblock {\em Ann. Inst. H. Poincar{\'e} Probab. Statist.}, {\bf
  39}(3):413--485, (2003).
}}


\bibitem{Holm08}
M.~Holmes.
\newblock Convergence of lattice trees to super-Brownian motion above the
  critical dimension.
  \newblock {\em Electron. J. Probab.}, {\bf  13}:671--755, (2008).


\bibitem{H08comp}
M.~Holmes.
\newblock Excited against the tide: A random walk with competing drifts.
\newblock Preprint, (2009).


\bibitem{H08rwre}
M.~Holmes and R.~Sun.
\newblock A monotonicity property for random walk in a partially random environment.
\newblock Preprint (2009).

\bibitem{IV08}
D.~Ioffe and Y.~Velenik.
\newblock Ballistic phase of self-interacting random walks.
\newblock In {\em Analysis and Stochastics of Growth Processes and Interface Models, P. Mörters et al. (eds)}
\newblock Oxford University Press,  55--79  (2008).

\bibitem{Kozm03}
G.~Kozma.
\newblock Excited random walk in three dimensions has positive speed.
\newblock Available on {\tt http://arxiv.org/abs/math.PR/0310305}, (2003).

\bibitem{Kozm05}
G.~Kozma.
\newblock Excited random walk in two dimensions has linear speed.
\newblock Available on {\tt http://arxiv.org/abs/math.PR/0512535}, (2005).

{\extend{\bibitem{NguYan93}
B.G. Nguyen and W-S. Yang.
\newblock Triangle condition for oriented percolation in high dimensions.
\newblock {\em Ann.\ Probab.}, {\bf 21}:1809--1844, (1993).
}}

{\extend{\bibitem{NguYan95}
B.G. Nguyen and W-S. Yang.
\newblock Gaussian limit for critical oriented percolation in high dimensions.
\newblock {\em J. Stat. Phys.}, {\bf 78}:841--876, (1995).
}}

\bibitem{Pema88}
R.~Pemantle.
\newblock Phase transition in reinforced random walk and {RWRE} on trees.
\newblock {\em Ann. Probab.}, {\bf 16}(3):1229--1241, (1988).

\bibitem{Pema07}
R.~Pemantle.
\newblock A survey of random processes with reinforcement.
\newblock {\em Probab. Surv.}, {\bf 4}:1--79 (electronic), (2007).

\bibitem{Roll02}
S.~Rolles.
\newblock {\em Random Walks in Stochastic Surroundings}.
\newblock PhD thesis, University of Amsterdam, (2002).

{\extend{\bibitem{Slad87}
G.~Slade.
\newblock The diffusion of self-avoiding random walk in high dimensions.
\newblock {\em Commun. Math. Phys.}, {\bf 110}:661--683, (1987).
}}

{\extend{\bibitem{Slad88}
G.~Slade.
\newblock Convergence of self-avoiding random walk to {B}rownian motion in high
  dimensions.
\newblock {\em J. Phys. A: Math. Gen.}, {\bf 21}:L417--L420, (1988).
}}

{\extend{\bibitem{Slad89}
G.~Slade.
\newblock The scaling limit of self-avoiding random walk in high dimensions.
\newblock {\em Ann. Probab.}, {\bf 17}:91--107, (1989).
}}

\bibitem{Zern05}
M.~Zerner.
\newblock Multi-excited random walks on integers.
\newblock {\em Probab. Th. Rel. Fields}, {\bf 133}:98--122, (2005).

\bibitem{Zern06}
M.~Zerner.
\newblock Recurrence and transience of excited random walks on {${\mathbb
  Z}^d$} and strips.
\newblock {\em Elect. Comm. in Probab.}, {\bf 11}:118--128, (2006).

\end{thebibliography}

\end{document}